
%
%
%
%
%
\newif\ifsect\newif\iffinal
\secttrue\finaltrue
\def\strutdepth{\dp\strutbox}

\def\lsimb#1{\vadjust{\vtop to0pt{\baselineskip\strutdepth\vss
    \llap{\ttt\string #1\ }\null}}}
\def\rsimb#1{\vadjust{\vtop to0pt{\baselineskip\strutdepth\vss
    \line{\kern\hsize\rlap{\ttt\ \string #1}}\null}}}
\def\ssect #1. {\bigbreak\indent{\bf #1.}\enspace\message{#1}}
\def\smallsect #1. #2\par{\bigbreak\noindent{\bf #1.}\enspace{\bf #2}\par
    \global\parano=#1\global\eqnumbo=1\global\thmno=1\global\esno=0
    \global\remno=0\global\remno=0\global\defno=0
    \nobreak\smallskip\nobreak\noindent\message{#2}}
\def\thm #1: #2{\medbreak\noindent{\bf #1:}\if(#2\thmp\else\thmn#2\fi}
\def\thmp #1) { (#1)\thmn{}}
\def\thmn#1#2\par{\enspace{\sl #1#2}\par
        \ifdim\lastskip<\medskipamount \removelastskip\penalty 55\medskip\fi}
\def\qedn{\thinspace\null\nobreak\hfill\hbox{\vbox{\kern-.2pt\hrule height.2pt
depth.2pt\kern-.2pt\kern-.2pt \hbox to2.5mm{\kern-.2pt\vrule width.4pt
\kern-.2pt\raise2.5mm\vbox to.2pt{}\lower0pt\vtop to.2pt{}\hfil\kern-.2pt
\vrule width.4pt\kern-.2pt}\kern-.2pt\kern-.2pt\hrule height.2pt depth.2pt
\kern-.2pt}}\par\medbreak}
\def\pf{\ifdim\lastskip<\smallskipamount \removelastskip\smallskip\fi
        \noindent{\sl Proof\/}:\enspace}
\def\itm#1{\par\indent\llap{\rm #1\enspace}\ignorespaces}

\def\bar#1{\overline{#1}}
\def\forclose#1{\hfil\llap{$#1$}\hfilneg}
\def\newforclose#1{
    \ifsect\xdef #1{(\number\parano.\number\eqnumbo)}\else
    \xdef #1{(\number\eqnumbo)}\fi
    \hfil\llap{$#1$}\hfilneg
    \global \advance \eqnumbo by 1
    \iffinal\else\rsimb#1\fi}
\def\newforclosea#1{
    \ifsect\xdef #1{{\rm(\the\Apptok.\number\eqnumbo)}}\else
    \xdef #1{(\number\eqnumbo)}\fi
    \hfil\llap{$#1$}\hfilneg
    \global \advance \eqnumbo by 1
    \iffinal\else\rsimb#1\fi}
\def\forevery#1#2$${\displaylines{\let\eqno=\forclose\let\neweqa=\newforclosea
        \let\neweq=\newforclose\hfilneg\rlap{$\qqquad\forall#1$}\hfil#2\cr}$$}
\newcount\parano
\newcount\eqnumbo
\newcount\thmno
\newcount\versiono
\versiono=0
\def\neweqt#1$${\xdef #1{(\number\parano.\number\eqnumbo)}
    \eqno #1$$
    \iffinal\else\rsimb#1\fi
    \global \advance \eqnumbo by 1}
\def\newthmt#1 #2: #3{\xdef #2{\number\parano.\number\thmno}
    \global \advance \thmno by 1
    \medbreak\noindent
    \iffinal\else\lsimb#2\fi
    {\bf #1 #2:}\if(#3\thmp\else\thmn#3\fi}
\def\neweqf#1$${\xdef #1{(\number\eqnumbo)}
    \eqno #1$$
    \iffinal\else\rlap{$\smash{\hbox{\hfilneg\string#1\hfilneg}}$}\fi
    \global \advance \eqnumbo by 1}
\def\newthmf#1 #2: #3{\xdef #2{\number\thmno}
    \global \advance \thmno by 1
    \medbreak\noindent
    \iffinal\else\llap{$\smash{\hbox{\hfilneg\string#1\hfilneg}}$}\fi
    {\bf #1 #2:}\if(#3\thmp\else\thmn#3\fi}
\def\inizia{\ifsect\let\neweq=\neweqt\else\let\neweq=\neweqf\fi
\ifsect\let\newthm=\newthmt\else\let\newthm=\newthmf\fi}
\def\bititolo{\empty}
\gdef\begin #1 #2\par{\xdef\titolo{#2}
\ifsect\let\neweq=\neweqt\else\let\neweq=\neweqf\fi
\ifsect\let\newthm=\newthmt\else\let\newthm=\newthmf\fi
\centerline{\titlefont\titolo}
\if\bititolo\empty\else\medskip\centerline{\titlefont\bititolo}
\xdef\titolo{\titolo\ \bititolo}\fi
\bigskip
\centerline{\bigfont \autore}
\if\istituto!\else\bigskip
\centerline{\istituto}
\centerline{\indirizzo}
\centerline{\email}\fi
\medskip
\centerline{#1~\anno}
\bigskip\bigskip
\ifsect\else\global\thmno=1\global\eqnumbo=1\fi}
\font\titlefont=cmssbx10 scaled \magstep1
\font\bigfont=cmr12

\font\bigbigslbf=cmbxsl10 scaled\magstep3
\font\bigbigfont=cmr10 scaled \magstep2
\font\ttt=cmtt10 at 10truept
\font\eightrm=cmr8
\font\eighttt=cmtt8

\let\sc=\smallcaps
\font\gto=eufm10
\font\eurb=eurb10
\def\gt #1{\hbox{\gto #1}}
\font\bbr=msbm10
\font\sbbr=msbm7
\font\ssbbr=msbm5
\def\ca #1{{\cal #1}}
\nopagenumbers
\binoppenalty=10000
\relpenalty=10000
\newfam\amsfam
\textfont\amsfam=\bbr \scriptfont\amsfam=\sbbr \scriptscriptfont\amsfam=\ssbbr
\let\de=\partial

\def\phe{\varphi}

\def\Hom{\mathop{\rm Hom}\nolimits}
\def\End{\mathop{\rm End}\nolimits}

\def\Im{\mathop{\rm Im}\nolimits}

\def\id{\mathop{\rm id}\nolimits}

\def\nub{\mathord{\hbox{\eurb\char"17}}}
\mathchardef\void="083F

\def\C{{\mathchoice{\hbox{\bbr C}}{\hbox{\bbr C}}{\hbox{\sbbr C}}
{\hbox{\sbbr C}}}}

\def\P{{\mathchoice{\hbox{\bbr P}}{\hbox{\bbr P}}{\hbox{\sbbr P}}
{\hbox{\sbbr P}}}}

\def\qqquad{\quad\qquad}

\newcount\notitle
\notitle=1
\headline={\ifodd\pageno\rhead\else\lhead\fi}
\def\rhead{\ifnum\pageno=\notitle\iffinal\hfill\else\hfill\tt Version
\the\versiono; \the\day/\the\month/\the\year\fi\else\hfill\eightrm\titolo\hfill
\folio\iffinal\else\eighttt /V.\the\versiono\fi\fi}
\def\lhead{\ifnum\pageno=\notitle\hfill\else\eightrm\folio\iffinal
\else{\eighttt /V.\the\versiono}\fi\hfill\autore\hfill
\fi}
\output={\plainoutput}
\newbox\bibliobox
\def\setref #1{\setbox\bibliobox=\hbox{[#1]\enspace}
    \parindent=\wd\bibliobox}
\def\biblap#1{\noindent\hang\rlap{[#1]\enspace}\indent\ignorespaces}
\def\art#1 #2: #3! #4! #5 #6 #7-#8 \par{\biblap{#1}#2: {\sl #3\/}.
    #4 {\bf #5} (#6)\if.#7\else, \hbox{#7--#8}\fi.\par\smallskip}
\def\book#1 #2: #3! #4 \par{\biblap{#1}#2: {\bf #3.} #4.\par\smallskip}
\def\coll#1 #2: #3! #4! #5 \par{\biblap{#1}#2: {\sl #3\/}. In {\bf #4,}
#5.\par\smallskip}
\def\pre#1 #2: #3! #4! #5 \par{\biblap{#1}#2: {\sl #3\/}. #4, #5.\par\smallskip}
%
%
\let\newthm=\newthmt
\let\neweq=\neweqt
\newcount\esno\newcount\defno\newcount\remno
\def\Def #1\par{\global \advance \defno by 1
    \medbreak
{\bf Definition \the\parano.\the\defno:}\enspace #1\par
\ifdim\lastskip<\medskipamount \removelastskip\penalty 55\medskip\fi}
\def\Rem #1\par{\global \advance \remno by 1
    \medbreak
{\bf Remark \the\parano.\the\remno:}\enspace #1\par
\ifdim\lastskip<\medskipamount \removelastskip\penalty 55\medskip\fi}
\def\Es #1\par{\global \advance \esno by 1
    \medbreak
{\sc Example \the\parano.\the\esno:}\enspace #1\par
\ifdim\lastskip<\medskipamount \removelastskip\penalty 55\medskip\fi}
\def\Bittersweet{} 
\def\raggedleft{\leftskip2cm plus1fill \spaceskip.3333em \xspaceskip.5em
\parindent=0pt\relax}
\def\NF #1\par{\medbreak\begingroup{}\raggedleft
#1\par\endgroup
\ifdim\lastskip<\medskipamount \removelastskip\penalty
55\medskip\fi} 
\def\Nota #1\par{\medbreak\begingroup\Bittersweet\raggedleft
#1\par\endgroup
\ifdim\lastskip<\medskipamount \removelastskip\penalty
55\medskip\fi}

\def\Downarrow{\Big\downarrow}
\def\longhookrightarrow{\lhook\joinrel\longrightarrow}
\def\mapright#1{\smash{\mathop{\longrightarrow}\limits^{#1}}}
\def\rmapdown#1{\Big\downarrow\rlap{$\vcenter{\hbox{$\scriptstyle#1$}}$}}

\def\lmapdown#1{\llap{$\vcenter{\hbox{$\scriptstyle#1$}}$}\Big\downarrow}

\def\autore{Marco Abate, Filippo Bracci, Francesca Tovena}
\def\istituto{!}
\def\indirizzo{}
\def\anno{2006}
\def\email{}
\versiono=5 \finaltrue
%
%
%
%
\pageno=0\notitle=0\null
\vfill
\centerline{\bigbigslbf Embeddings of submanifolds and normal bundles}
\bigskip
\centerline{\bigbigfont Marco Abate, Filippo Bracci, Francesca Tovena}
\bigskip
\centerline{\bigfont December 2006}
\bigskip\bigskip
\centerline{Marco Abate, e-mail: abate@dm.unipi.it}
\centerline{Dipartimento di Matematica, Universit\`a di
Pisa, Largo Pontecorvo 5, 56127 Pisa, Italy}
\smallskip
\centerline{Filippo Bracci, e-mail: fbracci@mat.uniroma2.it}
\centerline{Dipartimento di Matematica,
Universit\`a di Roma Tor Vergata, Via della Ricerca Scientifica, 00133 Roma, Italy}
\smallskip
\centerline{Francesca Tovena, e-mail: tovena@mat.uniroma2.it}
\centerline{Dipartimento di Matematica,
Universit\`a di Roma Tor Vergata, Via della Ricerca Scientifica, 00133 Roma, Italy}
\bigskip\bigskip
{\sc Abstract:} This paper is devoted to the study of the embeddings of a complex submanifold
$S$ inside a larger complex manifold $M$; in particular, we are interested in comparing the
embedding of $S$ in $M$ with the embedding of $S$  as the zero section in the total space of the
normal bundle~$N_S$ of~$S$ in $M$. We explicitely describe some cohomological classes allowing to
measure the difference between the two embeddings, in the spirit of the work by Grauert, Griffiths,
and Camacho-Movasati-Sad; we are also able to explain the geometrical meaning of the separate
vanishing of these classes. Our results holds for any codimension, but even for curves in a
surface we generalize previous results due to Laufert and Camacho-Movasati-Sad.
\vfill\vfill
\eject
\null\pageno=0
\vfill
\eject
\notitle=1
\begin {December} {Embeddings of submanifolds and normal bundles}

\smallsect 0. Introduction

This paper is devoted to the study of the embeddings of a complex submanifold
$S$ inside a larger complex manifold $M$; in particular, we are interested in comparing the
embedding of $S$ in $M$ with the embedding of $S$  as the zero section in the total space of the
normal bundle~$N_S$ of~$S$ in $M$. We explicitely describe some cohomological classes allowing to
measure the difference between the two embeddings, in the spirit of~[G, Gr, CMS]. Our hope is that
it may be a step towards a classification of foliations on $M$ transverse to a submanifold $S$;
see [CMS].

Our interest in this topic originated in our previous papers [ABT1, 2], where we studied index
theorems for holomorphic self-maps and foliations. We had a complex submanifold
$S$ of a complex manifold
$M$ and a holomorphic object
$\ca F$ (either a holomorphic self-map of $M$ fixing $S$ pointwise, or a possibly singular
holomorphic foliation of $M$); along the lines of the original Camacho-Sad index theorem [CS], we
wanted to recover Chern classes of the normal bundle $N_S$ of $S$ in $M$ by means of local
invariants associated to singular points of either $S$ or of the holomorphic object $\ca F$. It
turned out that to get index theorems of this kind one needs either hypotheses on the relative
position of $S$ and $\ca F$ (e.g., the holomorphic foliation should be tangent to $S$), or on the
embedding of $S$ into $M$: it should be close enough to the embedding of $S$ in~$N_S$ as zero
section.

We found two ways to express the geometrical conditions on the embedding we needed; either in
terms of the existence of local coordinates with suitable properties (in a way similar to what was
done in [CMS], a main source of inspiration for the present paper), or in a more intrinsic way, as
splittings of suitable exact sequence of sheaves, thus allowing us to rephrase the conditions in
terms of vanishing of cohomology classes. Furthermore, it turned out that we were actually working
only with the first two of a list of more and more stringent conditions on the embedding, and that
it might be interesting to study the whole list of conditions.

The first (well-known) condition on the embedding is the splitting condition. We say that $S$
{\sl splits} into~$M$ if the exact sequence
$$
O\longrightarrow TS\longrightarrow TM|_S \longrightarrow
N_S\longrightarrow O
$$
splits as sequence of vector bundles over~$S$, where $TS$ (respectively, $TM|_S$) is the
holomorphic tangent bundle of $S$ (respectively, of $M$ restricted to $S$). It turns out (see
Section~1) that
$S$ splits into~$M$ if and only if the exact sequence
$$
O\longrightarrow\ca I_S/\ca I_S^2\longrightarrow \ca O_M/\ca I_S^2\longrightarrow
\ca O_S=\ca O_M/\ca I_S\longrightarrow O
$$
splits as sequence of sheaves of rings, where $\ca O_M$ (respectively, $\ca O_S$) is the
structure sheaf of~$M$ (respectively, of $S$), and $\ca I_S$ is the ideal sheaf of~$S$.

Thus if $S$ splits we have a way to extend germs of holomorphic functions on~$S$ to germs of
holomorphic functions defined on $M$ up to the first order. It is then natural to say that
$S$ is {\sl $k$-splitting} into $M$ (for some~$k\ge 1$) if the exact sequence 
$$
O\longrightarrow\ca I_S/\ca I_S^{k+1}\longrightarrow\ca O_M/\ca I_S^{k+1}\longrightarrow
\ca O_S\longrightarrow O
$$
splits as sequence of sheaves of rings. If this happens, it turns out (see Section~3) that we can
introduce a structure of $\ca O_S$-module on $\ca I_S/\ca I_S^{h+1}$ for $2\le h\le k+1$ in such a
way that the sequences
$$
O\longrightarrow\ca I_S^h/\ca I_S^{h+1}\longrightarrow\ca I_S/\ca I_S^{h+1}\longrightarrow
\ca I_S/\ca I_S^h\longrightarrow O
$$
become exact sequences of $\ca O_S$-modules. If these sequences split, we say that $S$ is {\sl
$k$-comfortably embedded} in~$M$. (In [ABT1,~2] we introduced split, 2-split and 1-comfortably
embedded submanifolds only.)

We can characterize these conditions in terms of local coordinates. Indeed, in Section~2 we prove
that $S$ is $k$-splitting into $M$ if and only if there is an atlas $\gt
U=\{(U_\alpha,z_\alpha)\}$ of $M$ adapted to~$S$ (that is, such that $U_\alpha\cap
S\ne\void$ implies $U_\alpha\cap S=\{z_\alpha^1=\cdots=z_\alpha^m=0\}$, where $m$ is the
codimension of~$S$) such that
$$
\left.{\de^k z_\beta^p\over\de z_\alpha^{r_1}\cdots\de z_\alpha^{r_k}}\right|_S\equiv 0,
$$
for all $r_1,\ldots,r_k=1,\ldots, m$, all $p=m+1,\ldots,n=\dim M$, and all indeces $\alpha$,
$\beta$ such that $U_\alpha\cap U_\beta\cap S\ne\void$. Furthermore,
$S$ is $k$-comfortably embedded in~$M$ if and only if (Section~3) there is an atlas 
$\gt U=\{(U_\alpha,z_\alpha)\}$ of $M$ adapted to~$S$ such that
$$
\left.{\de^k z_\beta^p\over\de z_\alpha^{r_1}\cdots\de z_\alpha^{r_k}}\right|_S\equiv 0,
\quad\hbox{and}\quad
\left.{\de^{k+1} z_\beta^s\over\de z_\alpha^{r_1}\cdots\de z_\alpha^{r_{k+1}}}\right|_S\equiv 0,
$$
for all $r_1,\ldots,r_{k+1}$,~$s=1,\ldots,m$, all $p=m+1,\ldots,n=\dim M$, and all indeces
$\alpha$, $\beta$ such that $U_\alpha\cap U_\beta\cap S\ne\void$. In particular, we see that if
$S$ is $k$-splitting and $(k-1)$-comfortably embedded then we have an atlas~$\gt U$ such that the
changes of coordinates are of the form
$$
\cases{z_\beta^r=\sum_{s=1}^m (a_{\beta\alpha})^r_s(z''_\alpha)z_\alpha^s+R_{k+1},&for
$r=1,\ldots, m$,\cr
\noalign{\smallskip}
z_\beta^p=\phi^p_{\beta\alpha}(z''_\alpha)+R_{k+1},&for $p=m+1,\ldots,n$,
\cr}
$$
where $z''_\alpha=(z_\alpha^{m+1},\ldots,z_\alpha^n)$ are local coordinates on~$S$, and $R_{k+1}$
denotes a term belonging to~$\ca I_S^{k+1}$, that is vanishing of order at least~$k+1$ along~$S$.
We remark that in the total space of the normal bundle $N_S$ we can always find an atlas (the
natural one induced by any adapted atlas of $M$) with changes of coordinates of the form
$$
\cases{z_\beta^r=\sum_{s=1}^m (a_{\beta\alpha})^r_s(z''_\alpha)z_\alpha^s,&for
$r=1,\ldots, m$,\cr
\noalign{\smallskip}
z_\beta^p=\phi^p_{\beta\alpha}(z''_\alpha),&for $p=m+1,\ldots,n$;
\cr}
$$
therefore we may say that the embedding of a $k$-splitting and $(k-1)$-comfortably embedded
submanifold looks, up to order~$k$, like its embedding into the normal bundle.

This is the reason why we were led to study the classical
problem of comparing the embedding of a submanifold~$S$ in a
complex manifold~$M$ with the embedding of~$S$ in the normal bundle~$N_S$ as zero section. In
particular, one would like to know when these two embeddings are holomorphically equivalent, that
is when there is a neighbourhood~$U$ of $S$ in~$M$ biholomorphic to a neighbourhood of the zero
section in~$N_S$ (in some sense, such a $U$ would be a holomorphic tubular neighbourhood of~$S$
in~$M$). In most approaches to this problem (see, e.g., [G, An, Gr, Hi, CG, K1, K2, St, L, CMS, CM]
and references therein) the first step consists in showing that the two embeddings are
biholomorphic up to a finite order $k$ (and we shall say that $S$ is {\sl $k$-linearizable\/}) if
a suitable cohomology class vanishes. Then one gives geometrical conditions ensuring the
vanishing of all the involved cohomology groups, and that if the two embeddings are formally
isomorphic (that is, biholomorphic up to any finite order) then they are actually biholomorphic.

The main result of this paper is a direct proof (see Theorem~4.1) of the fact that a complex
submanifold is $k$-linearizable in this sense if and only if it is $k$-splitting and
$(k-1)$-comfortably embedded. Combining this with our results on the existence of suitable local
coordinates, we are then able to explicitely write two cohomology classes providing the
obstructions from $k$-linearizable to $(k+1)$-linearizable. It should be remarked that
most other authors (see, e.g., [G]) first gives this obstruction as a single cohomology class,
and then use a formal argument to split this class in two; our approach explains instead the
geometrical meaning of the independent vanishing of any of the two classes. Furthermore, our
results hold for any codimension, and not only in codimension one.

In Section 5 we exemplify our results in the case of a compact Riemann surface~$S$ embedded in a
complex surface~$M$. In particular, we are able to recover results originally proved in~[CMS]
under the slightly stronger assumption that $S$ is fibered embedded into~$M$ (which implies, in
particular, that $S$ is $k$-splitting for any~$k\ge 1$).

Finally, Section 6 is devoted to a slightly different characterization of 1-comfortably embedded
submanifolds. In [ABT1,~2] we showed that the 1-comfortably embedded condition can be used to
define holomorphic connections on suitable vector bundles; here we show that a possible
justification for this phenomenon is that 1-comfortably embedded is exactly equivalent to the
existence of an infinitesimal holomorphic connection on~$N_S$.

{\it Acknowledgements} We would like to thank Francesco Russo for
pointing out reference [MP], Jorge Vitorio Pereira for several useful conversations, and IMPA for
the warm hospitality during part of the preparation of this paper. This work has been partially
supported by INdAM and FIRB.

\smallsect 1. Holomorphic splitting

Let us begin by recalling some general terminology on exact sequences of sheaves.
We say that an exact sequence of sheaves (of abelian groups, rings, modules\dots)
$$
O\longrightarrow\ca R\mathop{\longrightarrow}\limits^\iota\ca S
\mathop{\longrightarrow}\limits^p \ca T\longrightarrow O
\neweq\eqseqr
$$
on a variety $S$ {\sl splits} if there is a morphism $\sigma\colon\ca T\to\ca S$ of sheaves
(of abelian groups, rings, modules\dots) such that $p\circ\sigma=\id$. Any such morphism is
called a {\sl splitting morphism.} It is easy to see that \eqseqr\ splits (as sequence of
sheaves of modules) if and only if there exists a {\sl left splitting morphism,} that is a
morphism of sheaves of modules~$\tau\colon\ca S\to\ca R$ such that
$\tau\circ\iota=\id$. Furthermore, for every splitting morphism $\sigma$ there exists a unique
left splitting morphism $\tau$ such that
$$
\iota\circ\tau+\sigma\circ p=\id.
\neweq\eqip
$$

Following Grothendieck and Atiyah (see [At]), one can give a cohomological characterization of
splitting for sequences of locally free $\ca O_S$-modules defined over a complex manifold~$S$. 

 Let $\ca E'$ and $\ca E''$ be two sheaves of locally free $\ca O_S$-modules over the
same complex manifold~$S$. An {\sl extension} of~$\ca E''$ by~$\ca E'$ is an exact sequence
of locally free $\ca O_S$-modules
$$
O\longrightarrow\ca E'\longrightarrow\ca E\longrightarrow\ca E''\longrightarrow O.
\neweq\eqseqE
$$
If $O\rightarrow\ca E'\rightarrow\tilde\ca E\rightarrow\ca E''\longrightarrow O$
is another extension of $\ca E''$ by~$\ca E'$, one says that the two extensions are {\sl
equivalent} if there is an isomorphism $\chi\colon\ca E\to\tilde\ca E$ of $\ca O_S$-modules
such that
$$

\matrix{O&\longrightarrow&\ca E'&\longrightarrow&\ca E&
\longrightarrow&\ca E''&\longrightarrow&O\cr
&&\Big\|&&\rmapdown{\chi}&&\Big\|&\cr
O&\longrightarrow&\ca E'&\longrightarrow&\tilde\ca E&
\longrightarrow&\ca E''&\longrightarrow&O\cr
}
$$
commutes.
 Almost by definition, an extension of $\ca E''$ by $\ca E'$ splits if and only if it is equivalent
to the trivial extension $O\to\ca E'\to\ca E'\oplus\ca E''\to\ca E''\to O$.

Applying the functor
$\Hom(\ca E'',\cdot)$ to the sequence \eqseqE\ one gets the exact sequence
$$
O\longrightarrow\Hom(\ca E'',\ca E')\longrightarrow\Hom(\ca E'',\ca E)\longrightarrow
\Hom(\ca E'',\ca E'')\longrightarrow O.
\neweq\eqseqEH
$$
Let $\delta\colon H^0\bigl(S,\Hom(\ca E'',\ca E'')\bigr)\to H^1\bigl(S,\Hom(\ca E'',\ca
E')\bigr)$ be the connecting homomorphism in the long exact cohomology sequence of~\eqseqEH.
Then one can associate to the exact sequence \eqseqE\ the cohomology class
$$
\delta(\id_{\ca E''})\in H^1\bigl(S,\Hom(\ca E'',\ca E')\bigr).
$$
This procedure gives a $1$-to-$1$ correspondance between the cohomology group
$H^1\bigl(S,\Hom(\ca E'',\ca E')\bigr)$ and isomorphism classes of extensions of~$\ca E''$ by~$\ca
E'$ (see [At, Proposition~1.2]):

\newthm Proposition \uotto: Let $S$ be a complex manifold, and $\ca E'$ and $\ca E''$ two locally
free $\ca O_S$-modules. Then two extensions of~$\ca E''$ by $\ca E'$ are equivalent if and only if
they correspond to the same cohomology class in~$H^1\bigl(S,\Hom(\ca E'',\ca E')\bigr)$. In
particular, the exact sequence $\eqseqE$ splits if and only if it corresponds to the zero
cohomology class.

Let us now introduce the sheaves we are interested in. Let $M$ be a complex manifold of
dimension $n$, and let $S$ be a reduced, globally irreducible subvariety of~$M$ of
codimension~$m\ge 1$. We denote: by $\ca O_M$ the sheaf of germs of holomorphic functions
on~$M$; by~$\ca I_S$ the subsheaf of~$\ca O_M$ of germs vanishing on~$S$; and by~$\ca O_S$
the quotient sheaf~$\ca O_M/\ca I_S$ of germs of holomorphic functions on~$S$.
Furthermore, let~$\ca T_M$ denote the sheaf of germs of holomorphic sections of the
holomorphic tangent bundle~$TM$ of~$M$, and $\Omega_M$ the sheaf of germs
of holomorphic 1-forms on~$M$. Finally, we shall denote by~$\ca T_{M,S}$ the sheaf of
germs of holomorphic sections along~$S$ of the restriction~$TM|_S$ of~$TM$ to~$S$, and by
$\Omega_{M,S}$ the sheaf of germs of holomorphic sections along~$S$ of~$T^*M|_S$. It is
easy to check that $\ca T_{M,S}=\ca T_M\otimes_{\ca O_M}\ca O_S$ and
$\Omega_{M,S}=\Omega_M\otimes_{\ca O_M}\ca O_S$. 

For $k\ge 1$ we shall denote by $f\mapsto[f]_k$ the canonical projection of $\ca O_M$ onto
$\ca O_M/\ca I_S^k$. The
{\sl cotangent sheaf}~$\Omega_S$ of~$S$ is defined by
$$
\Omega_S=\Omega_{M,S}\big/d_2(\ca I_S/\ca I_S^2),
$$
where $d_2\colon\ca O_M/\ca I_S^2\to\Omega_{M,S}$ is given by $d_2[f]_2=df\otimes[1]_1$.
In particular, we have the {\sl conormal sequence} of sheaves of~$\ca O_S$-modules associated
to~$S$:
$$
\ca I_S/\ca I_S^2\mathop{\longrightarrow}\limits^{d_2}\Omega_{M,S}
\mathop{\longrightarrow}\limits^{p}\Omega_S\longrightarrow O.
\neweq\eqsequno
$$
Applying the functor $\Hom_{\ca O_S}(\cdot,\ca O_S)$ to the conormal sequence we get the
{\sl normal sequence} of sheaves of $\ca O_S$-modules associated to~$S$:
$$
O\longrightarrow \ca T_S \mathop{\longrightarrow}\limits^{\iota}\ca T_{M,S}
\mathop{\longrightarrow}\limits^{p_2}\ca N_S,
\neweq\eqsqdue
$$
where $\ca T_S=\Hom_{\ca O_S}(\Omega_S,\ca O_S)$ is the {\sl tangent sheaf} of~$S$,
$\ca N_S=\Hom_{\ca O_S}(\ca I_S/\ca I_S^2,\ca O_S)$ is the {\sl normal sheaf} of~$S$, and
$p_2$ is the morphism dual to~$d_2$.

The first condition we shall consider on the embedding of the variety~$S$ inside~$M$ is:

\Def Let $S$ be a reduced, globally irreducible subvariety of a complex
manifold $M$. We say
that $S$ {\sl splits} into~$M$ if there exists a morphism of sheaves of
$\ca O_S$-modules
$\sigma\colon\Omega_S\to\Omega_{M,S}$ such that $p\circ\sigma=\id$, where
$p\colon\Omega_{M,S}\to\Omega_S$ is the canonical projection.

\Rem In the literature this notion has sometimes appeared under a different name; for instance,
Morrow and Rossi in [MR] say that the embedding $S\to M$ is direct.

It is not difficult to see that splitting subvarieties must be smooth:

\newthm Proposition \smooth: Let $S$ be a reduced, globally irreducible subvariety of a complex
manifold~$M$. Assume that $S$ splits in~$M$. Then $S$ is non-singular, and the morphism
$d_2\colon\ca I_S/\ca I_S^2\to\Omega_{M,S}$ is injective.

\pf Let us consider the sequence
$$
O\longrightarrow\ca K\longrightarrow\Omega_{M,S}
\mathop{\longrightarrow}\limits^{p}\Omega_S\longrightarrow O,
$$
where $\ca K=d_2(\ca I_S/\ca I_S^2)=\hbox{Ker}(p)$. This is a splitting exact sequence of $\ca
O_S$-modules; let $\sigma\colon\Omega_S\to\Omega_{M,S}$ be a splitting
morphism, and choose any $x\in S$. Then $\sigma\bigl((\Omega_S)_x\bigr)$ is a direct addend of the
projective module~$(\Omega_{M,S})_x$; therefore $\sigma\bigl((\Omega_S)_x\bigr)$ is itself
projective and thus, being~$\ca O_{S,x}$ Noetherian, it is $\ca O_{S,x}$-free. But $\sigma$
is an injective $\ca O_S$-morphism; therefore $(\Omega_S)_x$ itself is $\ca O_{S,x}$-free
for any $x\in S$.

Now, both $\Omega_{M,S}$ and $\Omega_S$  are coherent sheaves of
$\ca O_S$-modules, and $S$ is globally irreducible; therefore
$\Omega_S$ is locally $\ca O_S$-free of constant rank~$r$.
Checking the rank at a point in the regular part of~$S$ we see
that $r$ must be equal to the dimension of~$S$, and hence
[H, Theorem~I\negthinspace I.8.15] $S$ is non-singular,
and [H, Theorem~I\negthinspace I.8.17]~$d_2$ is
injective.\qedn

In particular, when $S$ splits into~$M$ the sequence
$$
O\longrightarrow\ca I_S/\ca I_S^2\mathop{\longrightarrow}\limits^{d_2}\Omega_{M,S}
\mathop{\longrightarrow}\limits^{p}\Omega_S\longrightarrow O
\neweq\eqsequnoprimo
$$
is a splitting exact sequence of $\ca O_S$-modules, and we also have a left splitting
morphism $\tau\colon\Omega_{M,S}\to\ca I_S/\ca I_S^2$. 

\Rem If $S$ splits into~$M$ then the normal sequence
$$
O\longrightarrow \ca T_S \mathop{\longrightarrow}\limits^{\iota}\ca T_{M,S}
\mathop{\longrightarrow}\limits^{p_2}\ca N_S\longrightarrow O
\neweq\eqseqdueprimo
$$
is a splitting exact sequence of $\ca O_S$-modules too: a splitting morphism is the dual
$\tau^*\colon\ca N_S\to\ca T_{M,S}$ of a left splitting morphism of~\eqsequnoprimo. Conversely,
if $S$ is a (reduced) locally complete intersection and \eqseqdueprimo\ is exact, then $S$ is
non-singular (see, e.g., [S]) and so if moreover \eqseqdueprimo\ splits then $S$ splits
into~$M$. There are examples of singular varieties for which \eqseqdueprimo\ is exact (see
again [S]); we do not know whether there are singular varieties for which \eqseqdueprimo\ is
exact and splits.

%

The aim of this section is to describe several equivalent characterizations of splitting
subvarieties. Most of them were already present in the literature; we collect them
here because they provide a template for the study of the more stringent conditions on the
embedding of~$S$ into~$M$ we shall study starting from the next section.

\Def Let $S$ be a reduced, globally irreducible subvariety of a
complex manifold $M$. For any~$k\ge h\ge 0$ let $\theta_{k,h}\colon\ca O_M/\ca I_S^{k+1}\to\ca O_M/\ca I_S^{h+1}$
be the canonical projection given by $\theta_{k,h}[f]_{k+1}=[f]_{h+1}$; when~$h=0$ we shall
write $\theta_k$ instead of~$\theta_{k,0}$. The {\sl $k$-th infinitesimal
neighbourhood} of $S$ in $M$ is the ringed space $S(k)=(S,\ca
O_M/\ca I_S^{k+1})$ together with the canonical inclusion of
ringed spaces $\iota_k\colon S=S(0)\to S(k)$ given by
$\iota_k=(\id_S,\theta_k)$. We shall also set $\ca O_{S(k)}=\ca O_M/\ca I_S^{k+1}$.

\Def A {\sl $k$-th order infinitesimal
retraction} is a morphism of ringed spaces $r\colon S(k)\to S$
such that $r\circ\iota_k=\id$. 
A $k$-th order infinitesimal retraction is given by a pair $r=(\id_S,\rho)$, where
$\rho\colon\ca O_S\to\ca O_{S(k)}$ is a  ring morphism  such that
$\theta_k\circ\rho=\id$. So, the existence of a $k$-th order infinitesimal retraction is equivalent
to the existence of a splitting morphism for the exact sequence of sheaves of rings
$$
O\longrightarrow\ca I_S/\ca I_S^{k+1}\longrightarrow\ca
O_{S(k)}\mathop{\longrightarrow}\limits^{\theta_k}\ca O_S\longrightarrow O.
\neweq\eqseqk
$$
Such a splitting morphism is called a {\sl $k$-th order lifting.}
More generally, for $k\ge h\ge 0$
we shall say that $S(k)$ {\sl retracts onto~$S(h)$} if there is a
morphism of ringed spaces $r\colon S(k)\to S(h)$ such that
$r\circ\iota_{h,k}=\id$, where
$\iota_{h,k}=(\id_S,\theta_{k,h})\colon S(h)\to S(k)$ is the
natural inclusion.

\Rem It is easy to see that
$\rho([1]_1)=[1]_{k+1}$ for any $k$-th order lifting
$\rho\colon\ca O_S\to\ca O_M/\ca I_S^{k+1}$. This is not
an automatic consequence of~$\rho$ being a morphism of sheaves of rings but can be
proved as follows: from $\theta_k\circ\rho=\id$ we get
$[1]_{k+1}-\rho([1]_1)\in\ca I_S/\ca I_S^{k+1}$, that is
$\rho([1]_1)=[1+h]_{k+1}$ for a suitable~$h\in\ca I_S$. Now
$$
[1+h]_{k+1}=\rho([1]_1)=\rho([1]_1)\rho([1]_1)=[(1+h)^2]_{k+1}=[1+2h+h^2]_{k+1},
$$
and so $[h+h^2]_{k+1}=O$. But $[1+h]_{k+1}$ is a unit in~$\ca
O_{S(k)}$; therefore $[h]_{k+1}=O$, and
$\rho([1]_1)=[1]_{k+1}$.

%
\Def Let $\ca O$, $\ca R$ be sheaves of rings, $\theta\colon\ca R\to\ca O$ a
morphism of sheaves of rings, and $\ca M$ a sheaf of $\ca O$-modules. A {\sl
$\theta$-derivation} of~$\ca R$ in~$\ca M$ is a morphism of sheaves of abelian groups
$D\colon\ca R\to\ca M$ such that
$$
D(r_1r_2)=\theta(r_1)\cdot D(r_2)+\theta(r_2)\cdot D(r_1)
$$
for any $r_1$, $r_2\in\ca R$. In other words, $D$ is a
derivation with respect to the $\ca R$-module structure induced via
restriction of scalars by~$\theta$.

We can now state a first list of properties equivalent to splitting,
(see [MP, Lemma~1.1] and [Ei, Proposition~16.12] for proofs)
including the existence of first order infinitesimal retractions:

\newthm Proposition \ucinque: Let $S$ be a reduced, globally irreducible
subvariety of a complex manifold~$M$. Then there is a $1$-to-$1$
correspondance among the following classes of morphisms:
\smallskip
\itm{(a)} morphisms $\sigma\colon\Omega_S\to\Omega_{M,S}$ of sheaves of
$\ca O_S$-modules such that $p\circ\sigma=\id$;
\itm{(b)} morphisms $\tau\colon\Omega_{M,S}\to\ca I_S/\ca I_S^2$ of sheaves of $\ca
O_S$-modules such that $\tau\circ d_2=\id$;
\itm{(c)} derivations $D\colon\ca O_M\to\ca I_S/\ca I_S^2$ such that
$D|_{\ca I_S}=\pi_2|_{\ca I_S}$;
\item{\rm(d)} morphisms $\tau_M\colon \Omega_M\to \ca I_S/\ca I_S^2$ of sheaves of $\ca
O_M$-modules such that $d_2\circ\tau_M= \pi$, where $\pi\colon\Omega_M\to\Omega_{M,S}$ is the
canonical projection;
\item{\rm (e)} $\theta_1$-derivations $\tilde\rho\colon\ca O_{S(1)}\to\ca I_S/\ca I_S^2$ such that
$\tilde\rho\circ i_1=\id$, where
$i_1\colon\ca I_S/\ca I_S^2\hookrightarrow\ca O_{S(1)}$ is the
canonical inclusion and $\theta_1\colon \ca O_{S(1)} \to \ca O_M/\ca I_S$ is the
canonical projection;
\itm{(f)} morphisms $\rho\colon\ca O_S\to\ca O_{S(1)}$ of sheaves of rings such that
$\theta_1\circ\rho=\id$.
\smallskip
\noindent In particular, 
$S$ splits into $M$ if and only if it admits a first order infinitesimal retraction.
Finally, if any (and hence all) of the classes {\rm (a)--(f)} is not empty, then it
is in $1$-to-$1$ correspondance with the following classes of morphisms
too:
\smallskip
\itm{(g)} morphisms $\tau^*\colon\ca N_S\to\ca T_{M,S}$ of sheaves of $\ca O_S$-modules such
that $p_2\circ\tau^*=\id$;
\itm{(h)} morphisms $\sigma^*\colon\ca T_{M,S}\to\ca T_S$ of sheaves of
$\ca O_S$-modules such that $\iota\circ\sigma^*=\id$.

We have already noticed that a splitting subvariety is necessarily non-singular; therefore
we can use differential geometric techniques to get another couple of
characterizations of splitting submanifolds.

\Def Let $S$ be a complex submanifold (not necessarily closed) of codimension $m\ge 1$ in an
$n$-dimensional complex manifold $M$, and let $(U_\alpha,z_\alpha)$ a chart of $M$. We shall
sistematically write
$z_\alpha=(z_\alpha^1,\ldots,z_\alpha^n)=(z'_\alpha,z''_\alpha)$, with
$z'_\alpha=(z_\alpha^1,\ldots,z_\alpha^m)$ and
$z''_\alpha=(z_\alpha^{m+1},\ldots,z_\alpha^n)$. We shall say that $(U_\alpha,z_\alpha)$ is
{\sl adapted} to~$S$ if either $U_\alpha\cap S=\void$ or $U_\alpha\cap
S=U_\alpha\cap\bar{S}=\{z_\alpha^1=\cdots=z_\alpha^m=0\}$. In particular, if
$(U_\alpha,z_\alpha)$ is adapted to~$S$ then $\{z_\alpha^1,\ldots,z_\alpha^m\}$ is a set of
generators of~$\ca I_{S,x}$ for all~$x\in U_\alpha\cap S$. An atlas $\gt
U=\{(U_\alpha,z_\alpha)\}$ of~$M$ is {\sl adapted} to~$S$ if all its charts are; then $\gt
U_S=\{(U_\alpha\cap S,z''_\alpha)\mid U_\alpha\cap S\ne\void\}$ is an atlas for~$S$. We shall say
that an atlas~$\{(U_\alpha,z_\alpha)\}$ adapted to~$S$ is {\sl projectable} if $z_\alpha\in
U_\alpha$ implies $(O',z''_\alpha)\in\ U_\alpha\cap S$ for any
$U_\alpha$ such that $U_\alpha\cap S\ne\void$. Clearly, every atlas adapted to $S$ can be refined
to a projectable adapted atlas. 

\Def Let $S$ be a complex submanifold (not necessarily closed) of codimension $m\ge 1$ in an
$n$-dimensional complex manifold $M$. The
{\sl normal bundle}~$N_S$ of~$S$ in~$M$ is the quotient bundle~$TM|_S/TS$; its dual is the
{\sl conormal bundle}~$N^*_S$. If
$(U_\alpha,z_\alpha)$ is a chart adapted to~$S$, for $r=1,\ldots,m$ we shall denote
by~$\de_{r,\alpha}$ the projection of~$\de/\de z_\alpha^r|_{U_\alpha\cap S}$ in~$N_S$, and
by~$\omega_\alpha^r$ the local section of~$N_S^*$ induced by~$dz_\alpha^r|_{U_\alpha\cap
S}$. Then $\{\de_{1,\alpha},\ldots,\de_{m,\alpha}\}$
and~$\{\omega_\alpha^1,\ldots,\omega_\alpha^m\}$ are local frames over~$U_\alpha\cap S$ for
$N_S$ and~$N_S^*$ respectively, dual to each other.

\Rem From now on, every chart and atlas we consider on~$M$ will be adapted to~$S$.
Furthermore, we shall use Einstein convention on the sum over repeated indices. Indices like
$j$, $h$, $k$ will run from~1 to~$n$; indices like $r$, $s$, $t$,~$u$, $v$ will run from~1
to~$m$; and indices like $p$, $q$ will run from~$m+1$ to~$n$.

\Rem If $(U_\alpha,z_\alpha)$ and $(U_\beta,z_\beta)$ are two adapted charts with
$U_\alpha\cap U_\beta\cap S\ne\void$, then it is easy to check that
$$
\left.{\de z_\beta^r\over\de z_\alpha^p}\right|_S\equiv O
\neweq\eqts
$$
for all $r=1,\ldots,m$ and $p=m+1,\ldots,n$.

\Def Let $\gt U=\{(U_\alpha,z_\alpha)\}$ be an adapted atlas for a complex submanifold~$S$
of codimension~$m\ge 1$ of a complex $n$-dimensional manifold~$M$. We say that $\gt U$ is a
{\sl splitting atlas} (see [ABT1, 2]) if
$$
\left.{\de z_\beta^p\over\de z_\alpha^r}\right|_S\equiv O
$$
for all $r=1,\ldots, m$, $p=m+1,\ldots,n$ and indices $\alpha$, $\beta$ such that
$U_\alpha\cap U_\beta\cap S\ne\void$. In other words, and recalling \eqts, the jacobian
matrices of the changes of coordinates become block-diagonal when restricted to~$S$.

\Def Let $\gt U=\{(U_\alpha,z_\alpha)\}$ be an atlas adapted to~$S$. If $\rho\colon\ca O_S\to\ca
O_{S(1)}$ is a first order lifting for~$S$, we say~$\gt U$ is
{\sl adapted} to~$\rho$ if
$$
\rho([f]_1)=[f]_2-\left[{\de f\over\de z^r_\alpha} z_\alpha^r\right]_2
\neweq\defgrho
$$
for all $f\in\ca O(U_\alpha)$ and all indices $\alpha$ such that $U_\alpha\cap S\ne\void$.

%
%
%
%

In [ABT2] we proved the following characterization of splitting submanifolds:

\newthm Proposition \udieci: Let $S$ be a complex submanifold of codimension~$m\ge 1$ of a
$n$-dimensional complex manifold~$M$. Then:
\smallskip
\itm{(i)}the cohomology class $\gt s\in H^1\bigl(S,\Hom(\Omega_S,\ca N^*_S)\bigr)$ associated to
the conormal exact sequence is represented\break\indent by the $1$-cocycle $\{\gt
s_{\beta\alpha}\}\in H^1\bigl(\gt U_S,\Hom(\Omega_S,\ca N^*_S)\bigr)$ given by
$$
\gt s_{\beta\alpha}=-\left.{\de z_\beta^r\over\de z_\alpha^s}\,{\de z_\alpha^p\over\de
z_\beta^r}\right|_S\omega_\alpha^s\otimes{\de\over\de z_\alpha^p}\in H^0(U_\alpha\cap
U_\beta\cap S,\ca N_S^*\otimes\ca T_S),
$$
\indent where $\gt U=\{(U_\alpha,z_\alpha)\}$ is an atlas adapted to $S$.
In particular, $S$ splits into $M$ if and only if $\gt s=O$;
\itm{(ii)}$S$ splits into $M$ if and only if there exists a splitting atlas
for $S$ in $M$;
\itm{(iii)}an atlas adapted to $S$ is splitting if and only if
it is adapted to a first order lifting;
\itm{(iv)}if $S$ splits into $M$, then for any first order lifting there exists an atlas
adapted to it.

\Rem Assume that $\gt U=\{(U_\alpha,z_\alpha)\}$ is a projectable atlas adapted to $S$.
Then if $f\in\ca O(U_\alpha)$ we can write
$$
f(z_\alpha)=f(O',z''_\alpha)+{\de f\over\de z_\alpha^r}(O',z''_\alpha)z_\alpha^r+R_2
=f(O',z''_\alpha)+{\de f\over\de z_\alpha^r}(z_\alpha)z_\alpha^r+R_2,
$$
where $R_2$ denotes an element of $\ca I_S^2(U_\alpha)$, possibly changing from one occurrence to
the next. From this formula it follows that
$\gt U$ is adapted to a first order lifting $\rho$ if and only if
$$
\rho([f]_1)=f(O',z''_\alpha)+R_2
$$
for every $f\in\ca O(U_\alpha)$. In other words, $\gt U$ is a splitting atlas if and only if we
can patch together the trivial local liftings $[f]_1\mapsto [f(O',z''_\alpha)]_2$ so to get a
global first order lifting.

\Rem Given a first order lifting $\rho\colon\ca O_S\to\ca O_M/\ca
I_S^2$ and an atlas $\gt
U=\{(U_\alpha,z_\alpha)\}$ adapted to~$S$, it is not
difficult to check that $\gt U$ is adapted to~$\rho$ if and only if for every
$(U_\alpha,z_\alpha)\in\gt U$ with
$U_\alpha\cap S\ne\void$ and every
$f\in\ca O_M|_{U_\alpha}$ one has
$$
\tilde\rho([f]_2)=\left[{\de f\over\de z_\alpha^r}z_\alpha^r\right]_2,
$$
where $\tilde\rho\colon\ca O_{S(2)}\to \ca I_S/\ca I_S^2$ is the $\theta_1$-derivation
associated to~$\rho$ by Proposition~\ucinque.

As mentioned in the introduction, one of the aims of our constructions will be the comparison of
the embedding of $S$ into $M$ with its embedding (as zero section) in the normal bundle. 
The first result of this kind is our last characterization of splitting submanifolds:

\newthm Proposition \usette: Let $S$ be a submanifold of a complex manifold~$M$. Then $S$
splits into~$M$ if and only if its first infinitesimal neighbourhood~$S(1)$ in~$M$ is isomorphic to its
first infinitesimal neighbourhood~$S_N(1)$ in~$N_S$, where we are identifying~$S$ with the
zero section of~$N_S$.

\pf The main observation here is that if $E$ is any vector bundle over~$S$, then $TE|_S$ is
canonically isomorphic to~$TS\oplus E$. When $E=N_S$ this implies that the projection
$TN_S|_S\to N_S$ on the second direct summand induces an isomorphism $N_{O_S}\to N_S$,
where $N_{O_S}$ is the normal bundle of~$S$ (or, more precisely, of the zero section
of~$N_S$) in~$N_S$; in particular, then,
$S$ always splits in~$N_S$ (see also Example~1.1 below). Furthermore, this isomorphism
induces an isomorphism between $\ca N_S^*$ and $\ca N_{O_S}^*$, and thus an isomorphism of
sheaves of
$\ca O_S$-modules $\chi\colon\ca I_{S,N_S}/\ca I_{S,N_S}^2\to\ca I_S/\ca I_S^2$, where
$\ca I_{S,N_S}$ is the ideal sheaf of $S$ in~$N_S$.

By definition, an isomorphism between $S_N(1)$ and $S(1)$ is given by an isomorphism of
sheaves of rings $\psi\colon\ca O_{N_S}/\ca I_{S,N_S}^2\to\ca O_M/\ca I_S^2$ such that
$\theta_1\circ\psi=\theta_1^N$, where $\theta_1^N\colon\ca O_{N_S}/\ca I_{S,N_S}^2\to\ca
O_S$
is the canonical projection.

If $S_N(1)$ and $S(1)$ are isomorphic, we can define a morphism of sheaves of rings
$\rho\colon\ca O_S\to\ca O_M/\ca I_S^2$ by setting $\rho=\psi\circ\rho^N$, where $\rho^N$
is the first order lifting induced by the splitting of~$S$ in~$N_S$ described above. Then it
is easy to see that~$\theta_1\circ\rho=\id$, and thus $S$ splits in~$M$ by Proposition~\ucinque.

Conversely, assume that $S$ splits in~$M$, and let $\rho\colon\ca O_S\to\ca O_M/\ca I_S^2$
be a first order lifting. Then we can define a morphism $\psi\colon\ca O_{N_S}/\ca
I_{S,N_S}^2\to\ca O_M/\ca I_S^2$ by setting
$$
\psi=\rho\circ\theta_1^N+i_1\circ\chi\circ \tilde\rho^N,
$$
where $\tilde\rho^N$ is the $\theta_1^N$-derivation associated to the first order
lifting~$\rho^N$ and $i_1\colon \ca I_S/\ca I_S^2\to \ca O_M/\ca I_S^2$ is the canonical
inclusion. Then it is not difficult to check that
$\psi$ is an isomorphism of sheaves of rings such that $\theta_1\circ\psi=\theta_1^N$, and
thus
$S_N(1)$ and $S(1)$ are isomorphic.\qedn

%
%
\Es A local holomorphic retract is always split in the ambient manifold (and thus it is
necessarily non-singular). Indeed, if $p\colon U\to S$ is a local holomorphic retraction,
then a first order lifting $\rho\colon\ca O_S\to\ca O_M/\ca I_S^2$ is given by
$\rho(f)=[f\circ p]_2$. In particular, the zero section of a vector bundle always splits, as
well as any slice $S\times\{x\}$ in a product $M=S\times X$ (with both $S$ and $X$
non-singular, of course).

%
\Es If $S$ is a Stein submanifold of a complex manifold~$M$ (e.g., if $S$ is an open Riemann
surface), then $S$ splits into~$M$. Indeed, we have $H^1(S,\ca T_S\otimes\ca N_S^*)=(O)$ by
Cartan's Theorem~B, and the assertion follows from Proposition~\udieci.(i). In particular, if $S$
is a {\it singular} curve in $M$ then the non-singular part of~$S$ always splits in~$M$.

%
%
%
\Es Let $S$ be a non-singular, compact, irreducible curve of genus~$g$ on a surface~$M$.
If $S\cdot S<4-4g$ then $S$ splits into~$M$. In fact, the Serre duality for
Riemann surfaces implies that
$$
H^1(S,\ca T_S\otimes\ca N_S^*)\cong H^0(S,\Omega_S\otimes\Omega_S\otimes\ca N_S),
$$
and the latter group vanishes because the line bundle $T^*S\otimes T^*S\otimes N_S$ has
negative degree by assumption. The bound $S\cdot S<4-4g$ is sharp: for instance, a
non-singular compact projective plane conic~$S$ has genus~$g=0$ and
self-intersection~$S\cdot S=4$, but it does not split in the projective plane (see
[VdV], [MR], [MP]).

\Es Let $M$ be an algebraic surface embedded in ${\bf P}^n$ and let $S$ be a section of
$M$ with an hyperplane $H$, with the property that there exists a point $P\not\in H$ 
belonging to each plane tangent to $M$ in points of $S$. Then $S$ splits in $M$. In [BM], the
authors show  a partial converse: 
if $S$ splits in $M$ and the natural morphism $H^0(S,\Omega_S)\otimes H^0({\bf P}^n,
\ca O_{{\bf P}^n})\to H^0(S,\Omega_S(1))$ is injective, then there exists
a point $P\not\in H$ 
belonging to each plane tangent to $M$ in points of $S$. 

\Es Let $S$ be a compact Riemann surface of genus $g>0$, and
$\phi\colon\pi_1(S)\to\hbox{Diff}_0(\C^n)$ be a representation of the fundamental group of~$S$
into the group of germs of biholomorphisms of~$\C^n$ fixing the origin; assume that all the
elements of the image of~$\phi$ are convergent on some polydisk~$\Delta\subseteq\C^n$ centered at
the origin. If
$\tilde S$ is the universal covering space of~$S$, we shall also identify $\pi_1(S)$ with the
group of the automorphisms of the covering. The {\sl suspension} $M$ of the representation~$\phi$
is by definition the quotient of~$\Delta\times\tilde S$ obtained identifying $(z,\tilde p)$ and 
 $(w,\tilde q)$ if and only if there exists $\gamma \in \pi_1(S)$ such that  $(w,\tilde q)
=(\rho(\gamma)(z),\gamma \cdot\tilde
p)$. Then $S$ embeds into~$M$ as the $0$-slice, that splits into $M$.

Other examples of splitting submanifolds are discussed in [ABT2].

\smallsect 2. $k$-splitting submanifolds

In the previous section we have seen that a complex submanifold $S$ of a complex manifold
$M$ splits into~$M$ if and only if the sequence
$$
O\longrightarrow\ca I_S/\ca I_S^2\mathop{\longrightarrow}\limits^{\iota_1}\ca O_M/\ca
I_S^2\mathop{\longrightarrow}\limits^{\theta_1}\ca O_M/\ca I_S\longrightarrow O
\neweq\eqlof
$$
splits as a sequence of sheaves of rings. This suggests a natural generalization:
%

\Def Let $S$ be a submanifold of a complex manifold~$M$, and $k\ge 1$. We shall say that $S$
{\sl $k$-splits} (or is {\sl $k$-splitting}) into~$M$ if there is an infinitesimal retraction
of~$S(k)$ onto~$S$, that is if there is a $k$-th order lifting~$\rho\colon\ca O_S\to\ca O_M/\ca
I_S^{k+1}$, or, in still other words, if the exact sequence
$$
O\longrightarrow\ca I_S/\ca I_S^{k+1}\longhookrightarrow\ca O_M/\ca
I_S^{k+1}\mathop{\longrightarrow}\limits^{\theta_k}\ca O_M/\ca I_S\longrightarrow O
\neweq\eqlofdue
$$
splits as sequence of sheaves of rings. 

\Rem In [Gr, p. 373] a $k$-splitting submanifold is called $k$-transversely foliated. 

The main result of this section is a characterization of $k$-splitting submanifolds along the
lines of Proposition~\udieci. To state it, we need the analogue of Definitions~1.7 and~1.8:

\Def Let $\gt U=\{(U_\alpha,z_\alpha)\}$ be an adapted atlas for a complex submanifold~$S$
of codimension~$m\ge 1$ of a complex $n$-dimensional manifold~$M$, and let $k\ge 1$. We say that
$\gt U$ is a {\sl $k$-splitting atlas} (see [ABT1, 2]) if
$$
{\de z_\beta^p\over\de z_\alpha^r}\in\ca I_S^k
\neweq\equkappa
$$
for all $r=1,\ldots, m$, $p=m+1,\ldots,n$ and indices $\alpha$, $\beta$ such that
$U_\alpha\cap U_\beta\cap S\ne\void$.

\Def We shall say that an atlas
$\{(U_\alpha,z_\alpha)\}$ adapted to~$S$ is {\sl adapted} to a $k$-th order
lifting~$\rho\colon\ca O_S\to\ca O_M/\ca I_S^{k+1}$ if
$$
\rho[f]_1=\sum_{l=0}^k (-1)^l {1\over l!}\left[{\de^l f\over
	\de z_\alpha^{r_1}\cdots\de z_\alpha^{r_l}}\, z_\alpha^{r_1}
\cdots z_\alpha^{r_l} \right]_{k+1},
\neweq\eqdue
$$
for every $f\in\ca O(U_\alpha)$ and all indices $\alpha$ such that
$U_\alpha\cap S\ne\void$.

Then: 



\newthm Theorem \pezzettino: Let $S$ be an $m$-codimensional submanifold of an
$n$-dimensional complex manifold~$M$. Then:
\smallskip
\item{\rm(i)} $S$ is $k$-splitting into~$M$ if and only if there exists a $k$-splitting atlas;
\item{\rm(ii)} an atlas adapted to $S$ is $k$-splitting if and only if it is adapted to a $k$-th
order lifting;
\item{\rm(iii)} a projectable atlas adapted to $S$ is $k$-splitting if and only if the local
$k$-th order liftings
$$
\rho_\alpha([f]_1)=f(O',z''_\alpha)+\ca I_S^{k+1}
\neweq\eqrhoalpha
$$
patch together to define a global $k$-th order lifting;
\item{\rm(iv)} if $S$ is $k$-splitting into~$M$ then every $k$-th order lifting admits an atlas
adapted to it.

\pf Let $\gt U=\{(U_\alpha,z_\alpha)\}$ be a projectable adapted atlas, and $f\in\ca O(U_\alpha)$;
first of all we would like to prove that 
$$
\rho_\alpha([f]_1)=\sum_{l=0}^k (-1)^l {1\over l!}\left[{\de^l f\over
	\de z_\alpha^{r_1}\cdots\de z_\alpha^{r_l}}\, z_\alpha^{r_1}
\cdots z_\alpha^{r_l} \right]_{k+1}.
\neweq\eqduerho
$$
Let us proceed by induction on $k$. For $k=1$ we have already proved this in Remark~1.7; so assume
that~\eqduerho\ holds for $k-1$. Then we can write
$$
\eqalign{
\rho_\alpha([f])_1&=[f]_{k+1}-\sum_{j=1}^k{1\over j!}\left[{\de^j f\over\de z_\alpha^{r_1}
\cdots\de z_\alpha^{r_j}}(O',z''_\alpha)z_\alpha^{r_1}\cdots z_\alpha^{r_j}\right]_{k+1}\cr
&=[f]_{k+1}-\sum_{j=1}^k{1\over j!}\sum_{h=0}^{k-j}(-1)^h{1\over h!}\left[{\de^{j+h} f\over\de
z_\alpha^{r_1}
\cdots\de z_\alpha^{r_{j+h}}}z_\alpha^{r_1}\cdots z_\alpha^{r_{j+h}}\right]_{k+1}\cr
&=[f]_{k+1}-\sum_{l=1}^k\left(\sum_{h=0}^{l-1}(-1)^h {l!\over h!(l-h)!}\right){1\over l!}
\left[{\de^l f\over\de z_\alpha^{r_1}
\cdots\de z_\alpha^{r_l}}z_\alpha^{r_1}\cdots z_\alpha^{r_l}\right]_{k+1}\cr
&=\sum_{l=0}^k (-1)^l {1\over l!}\left[{\de^l f\over
	\de z_\alpha^{r_1}\cdots\de z_\alpha^{r_l}}\, z_\alpha^{r_1}
\cdots z_\alpha^{r_l} \right]_{k+1},
\cr}
$$
as claimed. In particular, the right-hand side of \eqduerho\ is a ring morphism, and to get (ii)
it suffices to prove~(iii).

Let $\gt U=\{(U_\alpha,z_\alpha)\}$ be a projectable atlas adapted to~$S$, and assume there is
$0\le l\le k$ such that
%
%
$$
{\de z_\beta^p\over\de z_\alpha^r}\in\ca I_S^l
$$
for all $r=1,\ldots,m$ and $p=m+1,\ldots,n$, which is equivalent to
assuming that
$$
{\de^l z_\beta^p\over\de z_\alpha^{r_1}\cdots\de z_\alpha^{r_l}}\in\ca I_S
$$
for all $r_1,\ldots,r_l=1,\ldots,m$ and all $p=m+1,\ldots,n$. Then it easy to prove by
induction on~$l$ 
that we can write
$$
z^p_\alpha(z_\beta)=\phi_{\alpha\beta}^p
(z''_\beta)+h^p_{r_1\ldots r_{l+1}}(z_\beta)z_\beta^{r_1}\cdots z_\beta^{r_{l+1}}
\neweq\eqinpiu
$$
for suitable $\phi_{\alpha\beta}^p\in\ca O(U_\alpha\cap U_\beta\cap S)$
and $h^p_{r_a\ldots r_{l+1}}\in\ca O(U_\alpha\cap U_\beta)$, symmetric
in the lower indices;
clearly, $\phi_{\beta\alpha}\circ\phi_{\alpha\beta}=\id$.

To simplify the understanding of the subsequent computations, we shall explicitely use
the local chart $\phe_\alpha\colon U_\alpha\to\C^n$ associated to~$(U_\alpha,z_\alpha)$.
Now let $(U_\alpha,z_\alpha)$ and~$(U_\beta,z_\beta)\in\gt U$ 
such that $U_\alpha\cap U_\beta\cap S\ne\void$; we need to evaluate
$\rho_\beta([f]_1)-\rho_\alpha([f]_1)$. First of all we have
$$
\eqalign{
\rho_\beta([f]_1)-\rho_\alpha([f]_1)&=f\circ\phe_\beta^{-1}(O',z''_\beta)-
f\circ\phe_\alpha^{-1}(O',z''_\alpha)+\ca I_S^{k+1}\cr
&=f\circ\phe_\beta^{-1}(O',z''_\beta)-f\circ\phe_\beta^{-1}\Bigl(
\phe_\beta\circ\phe_\alpha^{-1}\bigl(O',(\phe_\alpha\circ\phe_\beta^{-1})''(z_\beta)
\bigr)\Bigr)+\ca I_S^{k+1}.\cr}
$$
%
Now, \eqinpiu\ yields
$$
\phe_\beta\circ\phe_\alpha^{-1}\bigl(O',(\phe_\alpha\circ\phe_\beta^{-1})''(z_\beta)
\bigr)=\left(O',\phi_{\beta\alpha}\bigl(\phi_{\alpha\beta}
(z''_\beta)+h_{r_1\ldots r_{l+1}}(z_\beta)z_\beta^{r_1}\cdots z_\beta^{r_{l+1}}
\bigr)\right)
$$
and
$$
\eqalign{
\phi_{\beta\alpha}\left(\phi_{\alpha\beta}
(z''_\beta)+h_{r_1\ldots r_{l+1}}(z_\beta)z_\beta^{r_1}\cdots z_\beta^{r_{l+1}}
\right)&=z''_\beta+{\de\phi_{\beta\alpha}\over\de z_\alpha^p}\bigl(\phi_{\alpha\beta}
(z''_\beta)\bigr)h^p_{r_1\ldots r_{l+1}}(z_\beta)z_\beta^{r_1}\cdots
z_\beta^{r_{l+1}}+R_{l+2}\cr &=z''_\beta+{\de z''_\beta\over\de z_\alpha^p}(O',z''_\alpha
)h^p_{r_1\ldots r_{l+1}}(z_\beta)z_\beta^{r_1}\cdots z_\beta^{r_{l+1}}+R_{l+2},
\cr}
$$
where, here and elsewhere, $R_j$ denotes a term with elements in~$\ca I_S^j$.
Therefore we get
$$
\eqalign{
\rho_\beta([f]_1)-\rho_\alpha([f]_1)&=-{\de z_\beta^q\over\de z_\alpha^p}
h^p_{r_1\ldots r_{l+1}}(z_\beta)z_\beta^{r_1}\cdots z_\beta^{r_{l+1}}{\de f\over\de
z_\beta^q}+R_{l+2}+\ca I_S^{k+1}\cr
&=-h^p_{r_1\ldots r_{l+1}}(z_\beta)z_\beta^{r_1}\cdots
z_\beta^{r_{l+1}}{\de f\over\de z_\alpha^p}+R_{l+2}+\ca I_S^{k+1}.
\cr}
$$

In particular, if $l=k$ we get $\rho_\alpha\equiv\rho_\beta$, and
thus if $\gt U$ is a $k$-splitting atlas we get a global $k$-th order lifting, proving one
direction in (i), (ii) and (iii). Conversely, if $l<k$ then $\rho_\alpha\not\equiv\rho_\beta$, and
thus we obtain the other direction in (ii) and (iii).

For later use, we explicitely remark that if $l=k-1$ then
$$
\rho_\beta-\rho_\alpha=-h^p_{r_1\ldots r_k}|_S\,{\de\over\de z^p_\alpha}\otimes
[z_\beta^{r_1}\cdots z_\beta^{r_k}]_{k+1}\in H^0(U_\alpha\cap U_\beta\cap S, \ca
T_S\otimes\ca I_S^k/
\ca I_S^{k+1}).
\neweq\eqell
$$
Furthermore, it is easy to see that
$$
h^p_{r_1\ldots r_k}|_S={1\over k!}\left. {\de z_\alpha^p\over\de z_\beta^{r_1}\cdots
\de z_\beta^{r_k}}\right|_S.
\neweq\eqecche
$$

Now, let us assume that we have a $k$-th order lifting $\rho\colon\ca O_S\to\ca
O_M/\ca I_S^{k+1}$; we claim that there exists an atlas adapted to~$\rho$. This will yield (iv)
and the missing direction in (i), completing the proof.

We shall argue by induction on~$k$. For $k=1$ the assertion follows from
Proposition~\udieci. Now let $k>1$. Then $\rho_1=\theta_{k,k-1}\circ\rho$ is a
$(k-1)$-th order lifting; let $\gt U=\{(U_\alpha,z_\alpha)\}$ be a (necessarily $(k-1)$-splitting)
projectable atlas adapted to~$\rho_1$. Define local $k$-th
order liftings~$\rho_\alpha$ as in~\eqrhoalpha, and set $\sigma_\alpha=\rho-\rho_\alpha$.
Now
$$
\theta_{k,k-1}\circ\sigma_\alpha=\rho_1-\theta_{k,k-1}\circ\rho_\alpha\equiv O,
$$
because the atlas is adapted to~$\rho_1$; therefore the image of~$\sigma_\alpha$ is
contained in~$\ca I_S^k/\ca I_S^{k+1}$. The latter is an~$\ca O_S$-module; we claim that
$\sigma_\alpha\colon\ca O_S|_{U_\alpha\cap S}\to\ca I_S^k/\ca I_S^{k+1}|_{U_\alpha\cap S}$
is a derivation. Indeed,
$$
\sigma_\alpha(fg)=\rho(f)\rho(g)-\rho_\alpha(f)\rho_\alpha(g)=\rho(f)(\rho-\rho_\alpha)(g)+
\rho_\alpha(g)(\rho-\rho_\alpha)(f)=f\cdot\sigma_\alpha(g)+g\cdot\sigma_\alpha(f),
$$
beacuse $\sigma_\alpha(f)$,~$\sigma_\alpha(g)\in\ca I_S^k/\ca I_S^{k+1}$, and
$uv=\theta_k(u)\cdot v$ for all $u\in\ca O_M/\ca I_S^{k+1}$ and $v\in\ca I_S^k/\ca
I_S^{k+1}$.

Hence we can find $(s_\alpha)^p_{r_1\ldots r_k}\in\ca O(U_\alpha\cap S)$, symmetric in
the lower indices, such that
$$
\sigma_\alpha=(s_\alpha)^p_{r_1\ldots r_k}{\de\over\de z_\alpha^p}\otimes
[z_\alpha^{r_1}\cdots z_\alpha^{r_k}]_{k+1}.
$$
Now, by construction $\sigma_\alpha-\sigma_\beta=\rho_\beta-\rho_\alpha$;
therefore \eqell\ yields
$$
h^p_{s_1\cdots s_k}
{\de z_\beta^{s_1}\over\de
z_\alpha^{r_1}}\cdots{\de z_\beta^{s_k}\over\de z_\alpha^{r_k}}+(s_\alpha)^p_{r_1\ldots
r_k}-
{\de z_\alpha^p\over\de z_\beta^q}(s_\beta)^q_{s_1\ldots s_k}{\de z_\beta^{s_1}\over\de
z_\alpha^{r_1}}\cdots{\de z_\beta^{s_k}\over\de z_\alpha^{r_k}}\in\ca I_S,
$$
and then
$$
h^p_{s_1\cdots
s_{k-1}r} {\de z_\beta^{s_1}\over\de
z_\alpha^{r_1}}\cdots{\de z_\beta^{s_{k-1}}\over\de z_\alpha^{r_{k-1}}}
+(s_\alpha)^p_{r_1\ldots r_k}{\de z_\alpha^{r_k}\over\de z_\beta^r}-
{\de z_\alpha^p\over\de z_\beta^q}(s_\beta)^q_{s_1\ldots s_{k-1}r}{\de z_\beta^{s_1}\over\de
z_\alpha^{r_1}}\cdots{\de z_\beta^{s_{k-1}}\over\de z_\alpha^{r_{k-1}}}\in\ca I_S.
\neweq\eqquasi
$$

Let us then consider the change of coordinates
$$
\cases{\hat z_\alpha^r=z_\alpha^r,\cr
\hat z_\alpha^p=z_\alpha^p+(s_\alpha)^p_{r_1\ldots r_k}(z''_\alpha)z_\alpha^{r_1}
\cdots z_\alpha^{r_k},\cr}
$$
defined in suitable open sets $\hat U_\alpha\subseteq U_\alpha$; we claim that $\{(\hat
U_\alpha,\hat z_\alpha)\}$ is the atlas we are looking for.
Indeed, we have
$$
\eqalign{
{\de\hat z^p_\alpha\over\de\hat z_\beta^r}&={\de\hat z^p_\alpha\over\de z_\beta^s}
{\de z^s_\beta\over\de\hat z_\beta^r}+{\de\hat z^p_\alpha\over\de z_\beta^q}
{\de z^q_\beta\over\de\hat z_\beta^r}={\de\hat z^p_\alpha\over\de z_\beta^r}+
{\de\hat z^p_\alpha\over\de z_\beta^q}
{\de z^q_\beta\over\de\hat z_\beta^r}\cr
&={\de z^p_\alpha\over\de z_\beta^r}+k\left[(s_\alpha)^p_{r_1\ldots r_k}z_\alpha^{r_1}
\cdots z_\alpha^{r_{k-1}}{\de z_\alpha^{r_k}\over\de z_\beta^r}-{\de z_\alpha^p\over
\de z_\beta^q}(s_\beta)^q_{s_1\ldots s_{k-1}r}z_\beta^{s_1}\cdots z_\beta^{s_{k-1}}
\right]+R_k.
\cr}
$$
Now, \eqinpiu\ with $l=k-1$ yields
$$
{\de z_\alpha^p\over\de z_\beta^r}=k h^p_{s_1\ldots s_{k-1} r}z_\beta^{s_1}\cdots
z_\beta^{s_{k-1}}+R_k,
$$
and so
$$
\eqalign{
{\de\hat z^p_\alpha\over\de\hat z_\beta^r}&=k\left[h^p_{s_1\ldots s_{k-1} r}z_\beta^{s_1}
\cdots z_\beta^{s_{k-1}}+(s_\alpha)^p_{r_1\ldots r_k}z_\alpha^{r_1}
\cdots z_\alpha^{r_{k-1}}{\de z_\alpha^{r_k}\over\de z_\beta^r}-{\de z_\alpha^p\over
\de z_\beta^q}(s_\beta)^q_{s_1\ldots s_{k-1}r}z_\beta^{s_1}\cdots z_\beta^{s_{k-1}}
\right]+R_k\cr
&=k\left[
h^p_{s_1\cdots
s_{k-1}r} {\de z_\beta^{s_1}\over\de
z_\alpha^{r_1}}\cdots{\de z_\beta^{s_{k-1}}\over\de z_\alpha^{r_{k-1}}}
+(s_\alpha)^p_{r_1\ldots r_k}{\de z_\alpha^{r_k}\over\de z_\beta^r}-
{\de z_\alpha^p\over\de z_\beta^q}(s_\beta)^q_{s_1\ldots s_{k-1}r}{\de z_\beta^{s_1}\over\de
z_\alpha^{r_1}}\cdots{\de z_\beta^{s_{k-1}}\over\de z_\alpha^{r_{k-1}}}
\right]z_\alpha^{r_1}\cdots z_\alpha^{r_{k-1}}+R_k\cr
&=R_k\in\ca I_S^k\cr}
$$
thanks to \eqquasi, where we used the fact that
$z^{r_k}_\alpha=(\de z^{r_k}_\alpha/\de z^r_\beta)z_\beta^r+R_2$.

Finally, we should check that $\{(\hat U_\alpha,\hat z_\alpha)\}$ is adapted to~$\rho$. But
indeed \eqell\ applied with $\hat z_\alpha$ instead of~$z_\beta$ yields
$$
f(O',\hat z''_\alpha)-f(O',z''_\alpha)=(s_\alpha)^p_{r_1\ldots r_k}
z_\alpha^{r_1}\cdots z_\alpha^{r_k}{\de f\over\de z_\alpha^p}+R_{k+1}=\sigma_\alpha(f)+
R_{k+1};
$$
hence
$$
\rho([f]_1)=\rho_\alpha([f]_1)+\sigma_\alpha([f]_1)=f(O',\hat z''_\alpha)+\ca I_S^{k+1},
$$
and the assertion follows from \eqduerho.
\qedn

\Rem
In particular, there is an infinitesimal retraction of~$S(k)$ onto~$S$ if and only if
there is an atlas $\{(U_\alpha,z_\alpha)\}$ adapted to~$S$ whose coordinates changes are of the
form
$$
\cases{z_\beta^r=(a_{\beta\alpha})^r_s(z_\alpha) z_\alpha^s&for $r=1,\ldots,m$,\cr
\noalign{\smallskip}
z_\beta^p=\phi^p_{\alpha\beta}(z''_\alpha)+R_{k+1}&for $p=m+1,\ldots,n$,\cr}
$$
which, roughly speaking, says that a neighbourhood of~$M$ is a fiber bundle over $S$ up to
order~$k$. The jets of the vector fields ${\de\over\de z^r_\alpha}$ in $\ca
T_M\otimes\ca O_{S(k)}$, for $r=1,\ldots,m$,  generate an infinitesimal foliation $\ca F_k $, i.e.,
an involutive submodule of~$\ca
T_M\otimes\ca O_{S(k)}$.


We explicitely compute the obstruction, predicted by [Gr, Prop. 1.6], for
passing from $(k-1)$-split to $k$-split:

\newthm Proposition \ksplit:
Let $S$ be an $m$-codimensional submanifold of an
$n$-dimensional complex manifold~$M$. 
Assume that $S$ is $(k-1)$-splitting in~$M$; let
$\rho_{k-1}\colon\ca O_S\to\ca O_M/\ca I_S^{k}$ be a $(k-1)$-th order
lifting, and $\{(U_\alpha,\phe_\alpha)\}$ a $(k-1)$-splitting atlas adapted to
$\rho_{k-1}$. Let $\gt g_k\in H^1\bigl(S,\Hom(\Omega_S,\ca
I_S^k/\ca I_S^{k+1})\bigr)$ be the cohomology class represented by the $1$-cocycle $\{(\gt
g_k)_{\beta\alpha}\}\in H^1\bigl(
\gt U_S,\Hom(\Omega_S,\ca I_S^k/\ca I_S^{k+1})\bigr)$ given by
$$
(\gt g_k)_{\beta\alpha}=-\left.{1\over k!}{\de^k z_\alpha^p\over\de z_\beta^{r_1}\cdots
\de z_\beta^{r_k}}
\right|_S\,{\de\over\de z^p_\alpha}\otimes [z_\beta^{r_1}\cdots z_\beta^{r_k}]_{k+1}\in
H^0(U_\alpha\cap U_\beta\cap S, \ca T_S\otimes\ca I_S^k/
\ca I_S^{k+1}).
$$
Then there exists a $k$-th order lifting $\rho_{k}\colon\ca O_S\to\ca O_M/\ca I_S^{k+1}$ such that
$\rho_{k-1} = \theta_{k-\rho_{k}, k}\circ \rho_{k}$ if and only if $\gt g_k =O$.

\pf One direction follows from the previous theorem, \eqell\ and \eqecche. Conversely, if $\gt g_k
=O$ up to shrinking the $U_\alpha$ we can find $(s_\alpha)^p_{r_1\ldots r_k}\in\ca O(U_\alpha\cap
S)$ such that setting
$$
\sigma_\alpha=(s_\alpha)^p_{r_1\ldots r_k}{\de\over\de z_\alpha^p}\otimes
[z_\alpha^{r_1}\cdots z_\alpha^{r_k}]_{k+1}.
$$
we get $(\gt g_k)_{\beta\alpha}=\sigma_\alpha-\sigma_\beta$. Then arguing as in the last part of
the proof of the previous theorem we find a $k$-splitting atlas, and we are done.\qedn

\smallsect 3. Comfortably embedded submanifolds

The sequence \eqlofdue\ is only one of the possible natural generalizations of
\eqlof. Another, apparently as natural, generalization is the sequence
$$
O\longrightarrow\ca I_S^2/\ca I_S^3\longhookrightarrow\ca O_M/\ca
I_S^3\mathop{\longrightarrow}\limits^{\theta_{2,1}}\ca O_M/\ca I_S^2\longrightarrow O;
\neweq\eqceszero
$$
the splitting (as sequence of sheaves of rings) of this exact sequence is equivalent to the
existence of an infinitesimal retraction of $S(2)$ onto~$S(1)$. Surprisingly enough,
this cannot ever happen:

\newthm Proposition \nosplit: Let $S$ be a reduced, globally irreducible subvariety of a
complex manifold~$M$, and take~$k>h\ge 1$. Assume there is an infinitesimal retraction
of~$S(k)$ onto~$S(h)$; then
$$
\lceil(k+1)/2\rceil< h+1
$$ 
(where
$\lceil x\rceil$ is the smallest integer greater than or equal to~$x$). In particular, there are no
infinitesimal retractions of~$S(k)$ onto~$S(1)$ for any $k\ge 2$.

\pf For any $1\le l\le h$ consider the following commutative diagram with exact rows and
columns
$$

\matrix{&&&&O&&O\cr
&&&&\Downarrow&&\Downarrow\cr
O&\longrightarrow&\ca
I_S^{h+1}/\ca I_S^{k+1}&\longhookrightarrow&\ca I_S^l/\ca I_S^{k+1}&
\mapright{\theta_{k,h}}&\ca I_S^l/\ca I_S^{h+1}&\longrightarrow&O\cr
&&\Big\|&&\Downarrow&&\Downarrow&\cr
O&\longrightarrow&\ca I_S^{h+1}/\ca I_S^{k+1}&
\longhookrightarrow&\ca O_M/\ca I_S^{k+1}&
\mapright{\theta_{k,h}}&\ca O_M/\ca I_S^{h+1}&\longrightarrow&O\cr
&&&&\lmapdown{\theta_{k,l-1}}&&\lmapdown{\theta_{h,l-1}}\cr
&&&&\ca O_M/\ca I_S^l&=\!=&\ca O_M/\ca I_S^l\cr
&&&&\Downarrow&&\Downarrow\cr
&&&&O&&O\cr
}.
$$
By assumption, we have a morphism of sheaves of rings $\rho\colon\ca O_M/\ca
I_S^{h+1}\to\ca O_M/\ca I_S^{k+1}$ such that~$\theta_{k,h}\circ\rho=\id$. Composing with
$\theta_{h,l-1}$ on the left we get
$$
\theta_{h,l-1}=\theta_{h,l-1}\circ\theta_{k,h}\circ\rho=\theta_{k,l-1}\circ\rho.
$$
This implies that $\rho(\ca I_S^l/\ca I_S^{h+1})\subseteq \ca I_S^l/\ca I_S^{k+1}$: indeed,
if $u\in\ca I_S^l/\ca I_S^{h+1}$ we have
$\theta_{k,l-1}\bigl(\rho(u)\bigr)=\theta_{h,l-1}(u)=O$, and hence $\rho(u)\in\ca I_S^l/\ca
I_S^{k+1}$.

Now, if $u\in\ca I_S^l/\ca I_S^{h+1}$ we have $u^r=O$ as soon as $r\ge(h+1)/l$. Therefore
if $r\ge(h+1)/l$ we have
$$
O=\rho(u^r)=\rho(u)^r\in \ca I_S^l/\ca I_S^{k+1}
$$
for all $u\in\ca I_S^l/\ca I_S^{h+1}$. But since $S$ is reduced,
we have $v^r=O$ in~$\ca I_S^l/\ca I_S^{k+1}$ if and only if
$v\in\ca I_S^p/\ca I_S^{k+1}$ with $p\ge(k+1)/r$. Therefore if
$\lceil (k+1)/r\rceil\ge h+1$  we have $\rho(\ca I_S^l/\ca
I_S^{h+1})\subseteq\ca I_S^{h+1}/\ca I_S^{k+1}$, and thus
$\theta_{k,h}\circ\rho|_{\ca I_S^l/\ca I_S^{h+1}}\equiv O$,
impossible; therefore, $\lceil (k+1)/r\rceil< h+1$.

Now, the largest value of $\lceil (k+1)/r\rceil$ is attained for the lowest value of~$r$;
and since $r\ge(h+1)/l$, the lowest value of~$r$ is~2, attained taking $l=h$. Therefore we
get $\lceil(k+1)/2\rceil<h+1$, as claimed.\qedn

%
%
%

The lesson suggested by the previous proof is that if one would like to study the splitting of
sequences of sheaves of rings like \eqceszero, it is important to take care of what happens in the
nilpotent part of the rings, that is in the sheaves $\ca I_S^h/\ca I_S^k$. 
We observe that the sheaf $\ca I_S^k/\ca I_S^{k+1}$ is isomorphic to the
symmetric power $\hbox{Sym}^k(\ca N_S^*)$ of the conormal sheaf, and thus it naturally is an $\ca
O_S$-module. The main new idea of this section is that when $S$ is
$k$-splitting then the sheaf $\ca I_S/\ca I_S^{k+1}$ too has a canonical structure of $\ca
O_S$-module:

\newthm Proposition \dunogen: Let $S$ be a complex submanifold of codimension~$m$ of a
complex manifold~$M$, and let~$\rho\colon\ca O_S\to\ca O_M/\ca
I_S^{k+1}$ be a $k$-th order lifting, with $k\ge 0$. Then for any $1\le h\le k+1$ the
lifting~$\rho$ induces a structure of locally $\ca O_S$-free module on~$\ca I_S/\ca
I_S^{h+1}$ so that  the sequence
$$
O\longrightarrow\ca I_S^h/\ca I_S^{h+1}\longhookrightarrow\ca
I_S/\ca I_S^{h+1}\mathop{\longrightarrow}\limits^{\theta_{h,h-1}}\ca I_S/\ca
I_S^h\longrightarrow O
\xdef \eqces{{$(\number\parano.\number\eqnumbo)_h$}}
    \eqno{(\number\parano.\number\eqnumbo)_h}
    \iffinal\else\rsimb\eqces\fi
    \global \advance \eqnumbo by 1
$$
 becomes an exact sequence of locally $\ca O_S$-free modules.

\pf We shall work by induction on~$k$. For $k=0$ there is nothing to prove; so let us assume
that the assertion holds for $k-1$. As we already remarked, for any $h\ge 1$ the sheaf $\ca
I_S^h/\ca I_S^{h+1}$ has a natural structure of locally free $\ca O_S$-module. The
$k$-th order lifting~$\rho$ induces a
$(k-1)$-order lifting $\rho_1=\theta_{k,k-1}\circ\rho$; therefore by induction for $1\le
h\le k$ we have a structure of locally free $\ca O_S$-module on~$\ca I_S/\ca I_S^{h+1}$ so
that all the \eqces\ become exact sequences of locally free~$\ca O_S$-modules.
Now, $\ca I_S/\ca I_S^{k+2}$ naturally is a $\ca O_M/\ca I_S^{k+1}$-module; we can then
endow
it with the $\ca O_S$-module structure obtained by restriction of the scalars via~$\rho$:
$$
v\cdot [h]_{k+2} = \rho(v)\cdot [h]_{k+2},
$$
for all $v\in\ca O_S$ and $h\in\ca I_S$, where in the right-hand
side we are using the $\ca O_M/\ca I_S^{k+1}$-module operation.
Since $\rho$ is a ring morphism and (by
Remark~1.4) $\rho[1]_1=[1]_{k+1}$, we get a well-defined structure of $\ca
O_S$-module on~$\ca I_S/\ca I_S^{k+2}$. We must verify that the inclusion
$\iota\colon\ca I_S^{k+1}/\ca I_S^{k+2}\hookrightarrow\ca I_S/\ca
I_S^{k+2}$ and the projection $\theta_{k+1,k}$ are $\ca
O_S$-module morphisms when $\ca I_S^{k+1}/\ca I_S^{k+2}$ has its
own $\ca O_S$-structure and $\ca I_S/\ca I_S^{k+1}$ has the $\ca
O_S$-structure induced by~$\rho_1$ by induction.

Given $v\in\ca O_S$, choose $f\in\ca O_M$ such that $v=[f]_1$, and $f^{\ca I}\in\ca I_S$
so that $\rho(v)=[f+f^{\ca I}]_{k+1}$.
Then for all~$g\in\ca I_S^{k+1}$ we have
$$
\iota(v\cdot[g]_{k+2})=\iota([fg]_{k+2})=[fg]_{k+2}=[(f+f^{\ca I})g]_{k+2}
=\rho(v)\cdot[g]_{k+2}=v\cdot\iota[g]_{k+2},
$$
and $\iota$ is an $\ca O_S$-morphism.

Analogously, if $g\in\ca I_S$ we have
$$
\theta_{k+1,k}(v\cdot[g]_{k+2})=\theta_{k+1,k}[fg+f^{\ca I}g]_{k+2}=[fg+f^{\ca
I}g]_{k+1}=\rho_1(v)\cdot[g]_{k+1} =v\cdot\theta_{k+1,k}[g]_{k+2},
$$
and $\theta_{k+1,k}$ is an $\ca O_S$-morphism.

Finally, since (by induction)
$\ca I_S/\ca I_S^{k+1}$ and $\ca I_S^{k+1}/\ca I_S^{k+2}$ are locally $\ca O_S$-free, $\ca
I_S/\ca I_S^{k+2}$ is locally $\ca O_S$-free too, and we are
done.\qedn

\Rem If $\{(U_\alpha,z_\alpha)\}$ is an atlas adapted to a $k$-th order lifting
$\rho\colon\ca
O_S\to\ca O_M/\ca I_S^{k+1}$, it is easy to see that
$\{[z_\alpha^r]_{h+1},[z_\alpha^{r_1}z_\alpha^{r_2}]_{h+1},\ldots, [z_\alpha^{r_1}\cdots
z_\alpha^{r_h}]_{h+1}\}$ is a local free set of generators of~$\ca I_S/\ca I_S^{h+1}$ over~$\ca
O_S$ for~$h=1,\ldots,k+1$.

We are thus led to the following generalization of the notion of comfortably embedded submanifolds
introduced in [ABT1, 2]:

\Def
Let $S$ be a (not necessarily closed) submanifold of a complex manifold~$M$, and let~$\rho\colon\ca
O_S\to\ca O_M/\ca I_S^{k+1}$ be a $k$-th order lifting, with $k\ge 1$. A {\sl comfortable splitting
sequence}~$\nub$ {\sl associated to}~$\rho$ is a $(k+1)$-uple
$\nub=(\nu_{0,1},\ldots,\nu_{k,k+1})$, where for $1\le h\le k+1$ each
$\nu_{h-1,h}\colon\ca I_S/\ca I_S^h\to\ca I_S/\ca I_S^{h+1}$ is a splitting $\ca
O_S$-morphism of the sequence \eqces\ with respect to the $\ca O_S$-module structures induced
by~$\rho$. A pair~$(\rho,\nub)$, where
$\rho$ is a $k$-th order lifting and $\nub$ is a comfortable splitting sequence
associated to~$\rho$, is called a {\sl $k$-comfortable pair} for $S$ in~$M$. We say that $S$
is {\sl $k$-comfortably embedded} in~$M$ with respect to~$\rho$ if it exists a $k$-comfortable
pair~$(\rho,\nub)$ for~$S$ in~$M$.

\Rem The choice of a $k$-comfortable pair $(\rho,\nub)$ fixes an isomorphism of $\ca O_S$-modules
$$
\ca O_M/\ca I_S^{k+1}\cong \ca O_S\oplus \ca I_S/\ca I_S^2 \oplus \ca I_S^2/\ca I_S^3\oplus
\cdots \oplus \ca I_S^k/\ca I_S^{k+1}.
$$

The computation of the cohomology class associated
to the exact sequence $(3.2)_h$ is not too difficult:

\newthm Proposition \cocycle: Let $S$ be a complex submanifold of codimension~$m$ of a complex
manifold~$M$, and let $\rho\colon\ca O_S\to\ca O_M/\ca I_S^{k+1}$ be a $k$-th order lifting.
Choose a projectable atlas~$\gt U=\{(U_\alpha,z_\alpha)\}$ adapted to~$\rho$. Then the $1$-cocycle
$\{\gt h_{\beta\alpha}^\rho\}\in H^1\bigl(\gt U_S,\Hom(\ca I_S/\ca I_S^{k+1},\ca
I_S^{k+1}/\ca I_S^{k+2})\bigr)$ given by
$$
\gt h_{\beta\alpha}^\rho([z_\alpha^{t_1}\cdots z_\alpha^{t_h}]_{k+1})=-\left.{1\over(k+1)!}
{\de z_\beta^{s_1}\over\de z_\alpha^{r_1}}\cdots{\de z_\beta^{s_{k+1}}\over\de
z_\alpha^{r_{k+1}}}{\de^{k+1}(z_\alpha^{t_1}\cdots z_\alpha^{t_h})\over\de z_\beta^{s_1}
\cdots\de z_\beta^{s_{k+1}}}\right|_S [z_\alpha^{r_1}\cdots z_\alpha^{r_{k+1}}]_{k+2}
\neweq\eqcocl
$$
for $1\le t_1,\ldots,t_h\le m$ and $1\le h\le k$, represents the
class~$\gt h^\rho\in H^1\bigl(S,\Hom(\ca I_S/\ca I_S^{k+1},\ca
I_S^{k+1}/\ca I_S^{k+2})\bigr)$ associated to the exact sequence
$$
O\to\ca I_S^{k+1}/\ca I_S^{k+2}\to\ca I_S/\ca I_S^{k+2}\to\ca I_S/\ca I_S^{k+1}\to O,
\neweq\eqextk
$$
where $\ca I_S/\ca I_S^{k+1}$ and $\ca I_S/\ca I_S^{k+2}$ have the $\ca O_S$-module structure
induced by $\rho$.

\pf We can define local splittings $\nu_\alpha
\colon \ca I_S/\ca I_S^{k+1}|_{U_\alpha}\to\ca I_S/\ca I_S^{k+2}|_{U_\alpha}$ by setting
$$
\nu_\alpha[z_\alpha^{t_1}\cdots z_\alpha^{t_h}]_{k+1}=[z_\alpha^{t_1}\cdots
z_\alpha^{t_h}]_{k+2}
$$
and extending by $\ca O_S$-linearity; since $\gt U$ is adapted to~$\rho$, Theorem~\pezzettino\
implies that each
$\nu_\alpha$ is a well-defined morphism of $\ca O_S$-modules.

Now, for any $f\in\ca O_M$ we can write
$$
f(z'_\beta,z''_\beta)=f(O',z''_\beta)+\sum_{j=1}^{k+1}{1\over j!}{\de^j f\over\de
z_\beta^{s_1}\cdots\de z_\beta^{s_j}}(O',z''_\beta) z_\beta^{s_1}\cdots
z_\beta^{s_j}+R_{k+2},
\neweq\eqnew
$$
where $R_{k+2}\in\ca I_S^{k+2}$. In particular, Theorem~\pezzettino.(iii) implies
$$
[h]_{k+2}=\sum_{j=1}^{k+1}{1\over j!}\,\rho\left(\left[{\de^j h\over\de
z_\beta^{s_1}\cdots\de z_\beta^{s_j}}\right]_1\right)\cdot [z_\beta^{s_1}\cdots
z_\beta^{s_j}]_{k+2},
$$
for all $h\in\ca I_S$, and
$$
[f]_{k+1}=\rho([f]_1)+\sum_{j=1}^k{1\over j!}\,\rho\left(\left[{\de^j f\over\de
z_\beta^{s_1}\cdots\de z_\beta^{s_j}}\right]_1\right) [z_\beta^{s_1}\cdots
z_\beta^{s_j}]_{k+1},
\neweq\eqestendi
$$
for all $f\in\ca O_M$.

Using these formulas we easily see that
$$
\eqalign{
\gt h^\rho_{\beta\alpha}([z_\alpha^{t_1}\cdots
z_\alpha^{t_h}]_{k+1})&=\nu_\beta([z_\alpha^{t_1}\cdots z_\alpha^{t_h}]_{k+1})-
\nu_\alpha([z_\alpha^{t_1}\cdots z_\alpha^{t_h}]_{k+1})
=\nu_\beta([z_\alpha^{t_1}\cdots z_\alpha^{t_h}]_{k+1})
-[z_\alpha^{t_1}\cdots z_\alpha^{t_h}]_{k+2}\cr
&=-{1\over (k+1)!}\left[{\de^{k+1} (z_\alpha^{t_1}\cdots z_\alpha^{t_h})\over\de
z_\beta^{s_1}\cdots\de z_\beta^{s_{k+1}}}\right]_1[z_\beta^{s_1}\cdots
z_\beta^{s_{k+1}}]_{k+2},\cr
&=-\left.{1\over(k+1)!}
{\de z_\beta^{s_1}\over\de z_\alpha^{r_1}}\cdots{\de z_\beta^{s_{k+1}}\over\de
z_\alpha^{r_{k+1}}}{\de^{k+1}(z_\alpha^{t_1}\cdots z_\alpha^{t_h})\over\de z_\beta^{s_1}
\cdots\de z_\beta^{s_{k+1}}}\right|_S [z_\alpha^{r_1}\cdots z_\alpha^{r_{k+1}}]_{k+2},
\cr}
$$
as claimed.\qedn

\newthm Corollary \kceobstr: Let $S$ be a $(k-1)$-comfortably embedded submanifold of a complex
manifold $M$, and let~$(\rho_1,\nub_1)$ be a $(k-1)$-comfortable pair. Assume that we have a $k$-th
order lifting~$\rho$ such that $\theta_{k,k-1}\circ\rho=\rho_1$. Then the sequence $\nub_1$
extends to a comfortable splitting sequence $\nub$ associated to~$\rho$ if and only if the
class~$\gt h^\rho\in H^1\bigl(S,\Hom(\ca I_S/\ca I_S^{k+1},
\ca
I_S^{k+1}/\ca I_S^{k+2})\bigr)$ vanishes.

We can characterize $k$-comfortably embedded submanifolds using adapted atlases.

\Def Let $S$ be a complex submanifold of codimension $m$ in a complex $n$-dimensional
manifold~$M$, and let $k\ge 1$. A {\sl $k$-comfortable atlas} is an atlas $\{(U_\alpha,z_\alpha)\}$
adapted to~$S$ such that 
$$
{\de z_\beta^p\over\de z_\alpha^r}\in\ca I_S^k\qquad\hbox{and}\qquad
{\de^2 z_\beta^r\over\de z_\alpha^{s_1}\de z_\alpha^{s_2}}\in\ca I_S^k
\neweq\equkappace
$$
for all $r$,~$s_1$,~$s_2=1,\ldots,m$, all $p=m+1,\ldots,n$ and all indices $\alpha$, $\beta$ such
that
$U_\alpha\cap U_\beta\cap S\ne\void$. In particular, a $k$-comfortable atlas is always
$k$-splitting.

\Def  Let $S$ be a $k$-comfortably embedded submanifold of codimension $m$ of a complex
manifold~$M$, and $(\rho,\nub)$ a $k$-comfortable pair for~$S$ in~$M$. We shall say
that an atlas~$\{(U_\alpha,z_\alpha)\}$ adapted to~$S$ is {\sl adapted}
to~$(\rho,\nub)$ if it is adapted to~$\rho$ and
$$
\nu_{h-1,h}([z_\alpha^r]_{h})=[z_\alpha^r]_{h+1}
$$
for all $1\le h\le k+1$, $1\le j\le h-1$ and $1\le r\le m$.

The following result, in particular, recovers the original definition of
1-comfortably embedded submanifold introduced in~[ABT]:

\newthm Theorem \dtre: Let $S$ be an $m$-codimensional submanifold of an
$n$-dimensional complex manifold~$M$. Then:
\smallskip
\item{\rm(i)} $S$ is $k$-comfortably embedded into~$M$ if and only if there exists a
$k$-comfortable atlas;
\item{\rm(ii)} an atlas adapted to $S$ is $k$-comfortable if and only if it is adapted to a
$k$-comfortable pair;
\item{\rm(iii)} if $S$ is $k$-comfortably embedded into~$M$ then every $k$-comfortable pair admits
a $k$-comfortable atlas\break\indent adapted to it.

\pf First of all, notice that \equkappace\ implies that
$$
{\de^l z_\beta^p\over\de z_\alpha^{r_1}\cdots\de z_\alpha^{r_l}}\in\ca
I_S^{k-l+1}\qquad\hbox{and}\qquad {\de^l z_\beta^r\over\de z_\alpha^{r_1}\cdots\de
z_\alpha^{r_l}}\in\ca I_S^{k-l+2}
\neweq\eqcompi
$$
for all $2\le l\le k+1$, $r$,~$r_1,\ldots,r_l=1,\ldots,m$, $p=m+1,\ldots,n$ and indices
$\alpha$, $\beta$ such that $U_\alpha\cap U_\beta\cap S\ne\void$.
In particular, if there exists a $k$-comfortable atlas then Theorem~\pezzettino,
Lemma~\cocycle\ and Proposition \uotto\ imply that $S$ is $k$-comfortably embedded, and that the
atlas is adapted to a $k$-comfortable pair.

To prove the rest of the theorem, we shall work by induction on~$k$. For $k=0$ there is nothing to
prove; so we assume that the statement holds for $k-1$, and that there exists a
$k$-comfortable pair $(\rho,\nub)$ for $S$. We must prove that it exists an atlas adapted
to~$(\rho,\nub)$, and that this atlas is necessarily $k$-comfortable.

Let $\rho_1=\theta_{k,k-1}\circ\rho$ and $\nub_1=(\nu_{0,1},\ldots,\nu_{k-1,k})$;
clearly,
$(\rho_1,\nub_1)$ is a $(k-1)$-comfortable pair for $S$ in~$M$. The induction
hypothesis then provides us with a $(k-1)$-comfortable atlas~$\gt U$ adapted
to~$(\rho_1,\nub_1)$; arguing as in the proof of Theorem~\pezzettino\ we can moreover modify this
atlas to get a new projectable $k$-splitting and $(k-1)$-comfortable atlas (still denoted by~$\gt
U$) adapted to~$\rho$ and to~$(\rho_1,\nub_1)$. We must now show how
to modify~$\gt U$ so to get an atlas adapted to~$(\rho,\nub)$, and to prove that such an atlas is
necessarily $k$-comfortable.

The first observation is that from $\de^2 z^r_\beta/\de z^s_\alpha\de z^t_\alpha\in\ca
I_S^{k-1}$ we get
$$
z^r_\beta=\left.{\de z_\beta^r\over\de z_\alpha^s}\right|_S z_\alpha^s + h^r_{s_1\ldots s_{k+1}}
z_\alpha^{s_1}\cdots z_\alpha^{s_{k+1}}+R_{k+2}
\neweq\eqabuno
$$
for suitable functions $h^r_{s_1\ldots s_{k+1}}\in\ca O(U_\alpha\cap U_\beta\cap S)$ symmetric
in the lower indices. Notice that $h^r_{s_1\ldots s_{k+1}}\equiv 0$ if and only if $\gt U$ is
$k$-comfortable. 

From \eqabuno\ we derive three identities that will be useful later:
$$
\eqalign{\displaystyle {\de z^r_\beta\over\de z_\alpha^s}&=\left.{\de z^r_\beta\over\de
z_\alpha^s}\right|_S+(k+1) h^r_{s_1\ldots s_kt}\,z_\alpha^{s_1}\cdots z_\alpha^{s_k}+R_{k+1},
\cr
\displaystyle{\de z^r_\beta\over\de z_\alpha^{s_1}\de z_\alpha^{s_2}}&=k(k+1) h^r_{r_1\ldots
r_{k-1}s_1s_2}\,z_\alpha^{r_1}\cdots z_\alpha^{r_{k-1}}+R_k,
\cr
z_\alpha^s&=\left.{\de z_\alpha^s\over\de z_\beta^r}\right|_S z_\beta^r -\left.
{\de z_\alpha^s\over\de z_\beta^r} h^r_{s_1\ldots
s_{k+1}}{\de z_\alpha^{s_1}\over\de z_\beta^{r_1}}\cdots
{\de z_\alpha^{s_{k+1}}\over\de z_\beta^{r_{k+1}}}\right|_S z_\beta^{r_1}\cdots
z_\beta^{r_{k+1}}+R_{k+2}.
\cr}
\neweq\eqmolte
$$
In particular it follows that
$$
\displaylines{
\quad[z_\alpha^{s_1}\cdots z_\alpha^{s_j}]_h\hfill\cr
\hfill=\cases{\displaystyle
\left.{\de z_\alpha^{s_1}\over\de z_\beta^{r_1}}\cdots
{\de z_\alpha^{s_j}\over\de z_\beta^{r_j}}\right|_S\cdot [z_\beta^{r_1}\cdots
z_\beta^{r_j}]_h&\hskip-4truecmif $1\le h\le k+2$, $1\le j\le h-1$, and $(j,h)\ne(1,k+2)$\cr
\noalign{\smallskip}
\displaystyle \left.{\de z_\alpha^{s_1}\over\de z_\beta^{r_1}}\right|_S\cdot
[z_\beta^{r_1}]_{k+2}-\left.
{\de z_\alpha^{s_1}\over\de z_\beta^r} h^r_{t_1\ldots
t_{k+1}}{\de z_\alpha^{t_1}\over\de z_\beta^{r_1}}\cdots
{\de z_\alpha^{t_{k+1}}\over\de z_\beta^{r_{k+1}}}\right|_S\cdot [z_\beta^{r_1}\cdots
z_\beta^{r_{k+1}}]_{k+2},&if $j=1$ and $h=k+2$.
\cr}}
$$

Now for every index $\alpha$ such that $U_\alpha\cap S\ne\void$ define
$\nu_{k,k+1;\alpha}\colon\ca I_S/\ca I_S^{k+1}|_{U_\alpha}\to \ca I_S/\ca
I_S^{k+2}|_{U_\alpha}$ by setting
$$
\nu_{k,k+1;\alpha}([z_\alpha^{s_1}\cdots z_{\alpha}^{s_j}]_{k+1})=[z_\alpha^{s_1}\cdots
z_{\alpha}^{s_j}]_{k+2}
$$
and then extending by $\ca O_S$-linearity. The previous computations imply that
$$
\nu_{k,k+1;\beta}([z_\alpha^{s_1}\cdots
z_{\alpha}^{s_j}]_{k+1})-\nu_{k,k+1;\alpha}([z_\alpha^{s_1}\cdots
z_{\alpha}^{s_j}]_{k+1})=
\cases{\displaystyle\left.
{\de z_\alpha^{s_1}\over\de z_\beta^r}\right|_S h^r_{r_1\ldots
r_{k+1}}\cdot [z_\alpha^{r_1}\cdots
z_\alpha^{r_{k+1}}]_{k+2}&if $j=1$,\cr
0& if $2\le j\le k$.\cr}
$$
In particular, $\gt U$ is adapted to~$(\rho,\nub)$ if and only if $h^r_{r_1\ldots r_{k+1}}\equiv
0$, if and only if $\gt U$ is $k$-comfortable.  

Now set
$\sigma_\alpha=\nu_{k,k+1}-\nu_{k,k+1;\alpha}$; since
$\theta_{k+1,k}\circ\sigma_\alpha=O$, it follows that $\Im\sigma_\alpha\subseteq\ca
I_S^{k+1}/\ca I_S^{k+2}$. In particular, there are $(c_\alpha)^s_{r_1\ldots r_{k+1}}\in\ca
O(U_\alpha\cap S)$, symmetric in the lower indices, such that
$$
\sigma_\alpha\bigl([z_\alpha^s]_{k+1}\bigr)=(c_\alpha)^s_{r_1\ldots r_{k+1}}\cdot
[z_\alpha^{r_1}\cdots z_\alpha^{r_{k+1}}]_{k+2}.
$$
Since $\sigma_\alpha-\sigma_\beta=\nu_{k,k+1;\beta}-\nu_{k,k+1;\alpha}$, we get
$$
\left.{\de z_\alpha^s\over\de z_\beta^r}\right|_S h^r_{r_1\ldots
r_{k+1}}\cdot [z_\alpha^{r_1}\cdots
z_\alpha^{r_{k+1}}]_{k+2}=\left[(c_\alpha)^s_{r_1\ldots r_{k+1}}-\left.{\de z_\alpha^s
\over\de z_\beta^r}(c_\beta)^r_{s_1\ldots s_{k+1}}{\de z_\beta^{s_1}\over\de z_\alpha^{r_1}}
\cdots{\de z_\beta^{s_{k+1}}\over\de z_\alpha^{r_{k+1}}}
\right|_S\right]\cdot [z_\alpha^{r_1}\cdots z_\alpha^{r_{k+1}}]_{k+2},
$$
that is
$$
h^r_{r_1\ldots r_{k+1}}+(c_\beta)^r_{t_1\ldots t_{k+1}}{\de z_\beta^{t_1}\over\de
z_\alpha^{r_1}}
\cdots{\de z_\beta^{t_{k+1}}\over\de z_\alpha^{r_{k+1}}}-{\de z_\beta^r\over\de z_\alpha^t}
(c_\alpha)^t_{r_1\ldots r_{k+1}}\in\ca I_S.
\neweq\equff
$$

We are finally ready to modify $\gt U$. We define new coordinates $\hat z_\alpha$ by setting
$$
\cases{\hat z^r_\alpha=z^r_\alpha+(c_\alpha)^r_{s_1\ldots s_{k+1}}(z''_\alpha)\,z_\alpha^{s_1}
\cdots z_\alpha^{s_{k+1}}&for $r=1,\ldots,m$,\cr
\hat z^p_\alpha=z^p_\alpha&for $p=m+1,\ldots,n$,
\cr}
$$
on suitable $\hat U_\alpha\subseteq U_\alpha$. We claim that $\hat{\gt U}=\{(\hat
U_\alpha,\hat z_\alpha)\}$ is as desired. First of all, it is easy to see that
$$
{\de\hat z_\beta^p\over\de\hat z^r_\alpha}={\de z_\beta^p\over\de z^r_\alpha}+R_k,
$$
and so $\hat{\gt U}$ is still $k$-splitting and adapted to~$\rho$. A quick
computation shows that
$$
{\de^2\hat z^r_\beta\over\de\hat z_\alpha^{s_1}\de\hat z_\alpha^{s_2}}=
{\de^2 z^t_\alpha\over\de\hat z_\alpha^{s_1}\de\hat z_\alpha^{s_2}}{\de\hat z_\beta^r\over
\de z_\alpha^t}+{\de z_\alpha^{t_1}\over\de\hat z_\alpha^{s_1}}{\de z_\alpha^{t_2}\over\de\hat
z_\alpha^{s_2}}{\de^2\hat z^r_\beta\over\de z_\alpha^{t_1}\de z_\alpha^{t_2}},
$$
and hence
$$
\displaylines{
{\de^2\hat z^r_\beta\over\de\hat z_\alpha^{s_1}\de\hat z_\alpha^{s_2}}
={\de^2 z^r_\beta\over\de z_\alpha^{s_1}\de z_\alpha^{s_2}}\hfill\cr
\hfill+k(k+1)\left[
(c_\beta)^r_{t_1\ldots t_{k+1}}{\de z^{t_1}_\beta\over\de z^{r_1}_\alpha}\cdots
{\de z^{t_{k-1}}_\beta\over\de z^{r_{t-1}}_\alpha}{\de z^{t_k}_\beta\over\de z^{s_1}_\alpha}
{\de z^{t_{k+1}}_\beta\over\de z^{s_2}_\alpha}-{\de z^r_\beta\over\de z_\alpha^t}
(c_\alpha)^t_{r_1\ldots r_{k-1}s_1s_2}
\right]z_\alpha^{r_1}\cdots z_\alpha^{r_{k-1}}+R_k\in\ca I_S^k\quad
\cr}
$$
as desired, thanks to \eqmolte\ and \equff.

Finally, it is easy to check that $\hat{\gt U}$ is adapted to~$(\rho,\nub)$.
Indeed, if we define $\hat\sigma_\alpha$ by using $\hat{\gt U}$ instead of~$\gt U$, the
previous calculations can be used to show that
$\sigma_\alpha-\hat\sigma_\alpha=\sigma_\alpha$, and thus $\hat\sigma_\alpha=O$, which means
exactly (recalling the induction hypothesis) that
$\hat{\gt U}$ is adapted to~$(\rho,\nub)$. \qedn

\Rem In other words, $S$ is $k$-comfortably embedded into~$M$ if and only if there is
an atlas~$\{(U_\alpha,z_\alpha)\}$ adapted to~$S$ whose changes of coordinates
are of the form
$$
\cases{z_\beta^r=(a_{\beta\alpha})^r_s(z''_\alpha) z_\alpha^s+R_{k+2}&for $r=1,\ldots,m$,\cr
\noalign{\smallskip}
z_\beta^p=\phi^p_{\alpha\beta}(z''_\alpha)+R_{k+1}&for $p=m+1,\ldots,n$.\cr}
\neweq\eqkce
$$
%

As a corollary of the previous theorem, we are able to characterize the obstruction for
passing from $(k-1)$-comfortably embedded to $k$-comfortably embedded:

\newthm Corollary \obsce: Let $S$ be an $m$-codimensional $k$ split submanifold of an
$n$-dimensional complex manifold~$M$ and assume that  there exists a $k$-th order lifting
$\rho\colon\ca O_S\to\ca O_M/\ca I_S^{k+1}$ such that $S$ is $(k-1)$-comfortably embedded in~$M$
with respect to~$\rho_1=\theta_{k,k-1}\circ\rho$. Fix a $(k-1)$-comfortable pair
$(\rho_1,\nub_1)$, and let
$\gt U=\{(U_\alpha,z_\alpha)\}$ be a projectable atlas adapted to~$\rho$ and~$(\rho_1,\nub_1)$.
 Then
the cohomology class $\gt h^\rho$ associated to the exact sequence $\eqextk$ is represented by the
$1$-cocycle $\{\tilde{\gt h}^\rho_{\beta\alpha}\}\in H^1(\gt U_S,\ca N_S\otimes\ca I_S^{k+1}/\ca
I_S^{k+2})$ given by
$$
\tilde{\gt h}^\rho_{\beta\alpha}=-\left.{1\over(k+1)!}{\de z_\beta^{s_1}\over\de
z_\alpha^{r_1}}\cdots {\de z_\beta^{s_{k+1}}\over\de z_\alpha^{r_{k+1}}}{\de^{k+1}
z_\alpha^t\over\de z_\beta^{s_1}\cdots
\de z_\beta^{s_{k+1}}}\right|_S\,\de_{t,\alpha}\otimes[z_\alpha^{r_1}\cdots
z_\alpha^{r_{k+1}}]_{k+2}.
\neweq\eqcoclb
$$
Thus $S$ is $k$-comfortably embedded (with respect to~$\rho$) if and only if $\gt
h^\rho=O$ in $H^1\bigl(S,\ca N_S\otimes\hbox{\rm Sym}^{k+1}(\ca N_S^*)\bigr)$. 

\pf The $(k-1)$-comfortable pair $(\rho_1,\nub_1)$ induces a canonical splitting 
$$
\ca I_S/\ca I_S^{k+1}\cong \bigoplus_{h=1}^k \ca I_S^h/\ca I_S^{h+1};
$$
therefore the class $\gt h^\rho$ associated to the sequence \eqextk\ and computed in
Proposition~\cocycle\ lives in
$$
H^1\bigl(S,\Hom(\ca I_S/\ca I_S^{k+1},\ca I_S^{k+1}/\ca I_S^{k+2})\bigr)\cong
\bigoplus_{h=1}^k H^1\bigl(S,\Hom(\ca I_S^h/\ca I_S^{h+1},\ca I_S^{k+1}/\ca I_S^{k+2})\bigr).
$$
The expression of $\gt h^\rho$ given in \eqcocl\ clearly reflects this decomposition. 
Now, \eqcompi\ implies that
$$
{\de^l z_\alpha^r\over\de z_\beta^{r_1}\cdots\de
z_\beta^{r_l}}\in\ca I_S
$$
for all $2\le l\le k$. Therefore \eqcocl\ shows that the only non-zero component of $\gt h^\rho$ 
is the one contained in  
$$
H^1\bigl(S,\Hom(\ca I_S/\ca I_S^2,\ca I_S^{k+1}/\ca I_S^{k+2})\bigr)\cong
H^1\bigl(S,\ca N_S\otimes\ca I_S^{k+1}/\ca I_S^{k+2}\bigr)\cong 
H^1\bigl(S,\ca N_S\otimes\hbox{Sym}^{k+1}(\ca N_S^*)\bigr),
$$
and its expression is given by \eqcoclb.\qedn 

Recalling Proposition~\ksplit, we then see that the obstruction for passing from $(k-1)$-split to
$k$-split lives in $H^1\bigl(S,\ca T_S\otimes\hbox{Sym}^k(\ca N_S^*)\bigr)$, while the
obstruction for passing from $(k-1)$-comfortably embedded to $k$-comfortably embedded lives in
$H^1\bigl(S,\ca N_S\otimes\hbox{Sym}^{k+1}(\ca N_S^*)\bigr)$. Now, a vanishing theorem due to
Grauert ([G, Hilfssatz~1]; see also [CM]) says that if $N_S$ is negative in the sense of Grauert
(that is, the zero section of~$N_S$ can be blown down to a point) then these groups vanish for $k$
large enough. We thus obtain the following

\newthm Corollary \tutform: Let $S$ be an $m$-codimensional compact complex submanifold of an
$n$-dimensional manifold~$M$, and assume that $N_S$ is negative in the sense of Grauert. Then there
exists a~$k_0\ge 1$ such that if~$S$ is
$k_0$-splitting (respectively, $k_0$-comfortably embedded) in~$M$ then it is $k$-splitting
(respectively, $k$-comfortably embedded) for all $k\ge k_0$.

A similar result can also be obtained assuming instead that $N_S$ is positive in a suitable sense;
see~[Gr, CG, K1, K2, St].

\Rem
At present we do not know whether a submanifold which is $k$-comfortably embedded with respect to
a given $k$-th order lifting is $k$-comfortably embedded with respect to any $k$-th order
lifting.

We end this section with some
examples of $k$-split and $k$-comfortably embedded submanifolds. We refer to section 5 for a more
detailed study of $k$-split and $k$-comfortably embedded curves in a surface.

\Es The zero section of a vector bundle is always $k$-split and $k$-comfortably embedded  in the
total space of the bundle for any $k\ge1$: indeed, any atlas trivializing the bundle satisfies
$$
{\de z_\beta^p\over\de z_\alpha^r}\equiv {\de^2 z_\beta^r\over\de z_\alpha^s\de z_\alpha^t}
\equiv 0.
$$

\Es A local holomorphic retract is always $k$-split in the
ambient manifold. Indeed, if $p\colon U\to S$ is a
local holomorphic retraction, then a $k$-th order lifting $\rho\colon\ca O_S\to\ca O_M/\ca
I_S^{k+1}$ is given by
$\rho(f)=[f\circ p]_{k+1}$.

\Es If $S$ is a Stein submanifold of a complex manifold~$M$ (e.g., if $S$ is an open Riemann
surface), then $S$ is $k$-split and $k$-comfortably embedded into~$M$ for any~$k\ge 1$.
Indeed, by Cartan's Theorem~B the first cohomology group of~$S$ with
coefficients in any coherent sheaf vanishes, and the assertion follows from
Propositions~\uotto, \ksplit\ and~\cocycle. In particular, if
$S$ is a {\it singular} curve in $M$ then the non-singular part of~$S$ is always comfortably
embedded in~$M$.

\Es Let $\tilde M$ be the blow-up of a submanifold~$X$ in a complex manifold~$M$. Then the
exceptional divisor~$E\subset\tilde M$ is $k$-split and $k$-comfortably embedded in~$\tilde M$
for any~$k\ge 1$: indeed, it is easy to  check that the atlas of~$\tilde M$ induced by an atlas
of~$M$ adapted to~$X$ is a $k$-comfortable atlas for any $k\ge 1$.

\smallsect 4. Embeddings in the normal bundle and $k$-linearizable submanifolds

Proposition~\usette\ suggests a third way of generalizing the notion of splitting submanifold:

\Def Let $S$ be a complex submanifold of a complex manifold~$M$. We shall say that $S$ is
{\sl $k$-linearizable} if its $k$-th infinitesimal neighbourhood~$S(k)$ in~$M$ is isomorphic
to its $k$-th infinitesimal neighbourhood~$S_N(k)$ in~$N_S$, where we are identifying~$S$
with the zero section of~$N_S$.

We have seen that $S$ splits into~$M$ if and only if it is 1-linearizable
(Proposition~\usette). In general, however, $k$-split does not imply $k$-linearizable
(while the converse hold). The missing link is provided by the notion of $(k-1)$-comfortably
embedded:

\newthm Theorem \missinglink: Let $S$ be a complex submanifold of a complex manifold~$M$,
and $k\ge 2$. Then $S$ is $k$-linearizable if and only if it is $k$-split and
$(k-1)$-comfortably embedded (with respect to the $(k-1)$-th order lifting induced by the
$k$-splitting).

\pf We shall denote by $\ca I_{S,N}$ the ideal sheaf of~$S$ in~$N_S$, and by
$\theta^N_{h,k}\colon\ca O_{N_S}/\ca I_{S,N}^{h+1}\to\ca O_{N_S}/\ca I_{S,N}^{k+1}$ the
canonical projections. Notice that $S$ is~$k$-split and $k$-comfortably embedded
in~$N_S$ for any~$k\ge 1$, by, for instance, Example~3.1;
we shall denote by $\rho_k^N\colon\ca O_S\to\ca O_{N_S}/\ca I_{S,N}^{k+1}$ and
$\nu_{k-1,k}^N\colon\ca I_{S,N}/\ca I_{S,N}^k\to\ca I_{S,N}/\ca I_{S,N}^{k+1}$
the corresponding morphisms.


We shall work by induction on~$k$. We have already seen that $S_N(1)\cong S(1)$ implies that
$S$ is 1-split. Suppose now that $S_N(k-1)\cong S(k-1)$ implies that $S$ is
$(k-1)$-split and $(k-2)$-comfortably embedded, and assume that
$S_N(k)\cong S(k)$. Let
$\psi\colon\ca O_{N_S}/\ca I^{k+1}_{S,N}\to\ca O_M/\ca I_S^{k+1}$
be a ring isomorphism such that~$\theta_k\circ\psi=\theta_k^N$.
The gist of the proof is contained in the following commutative
diagrams:

$$

\displaylines{\matrix{
O&\longrightarrow&\ca
I_{S,N}/\ca I_{S,N}^{k+1}&\longhookrightarrow&\ca O_{N_S}/\ca I_{S,N}^{k+1}&
{\textstyle\smash{\mathop{\longrightarrow}\limits^{\theta_k^N}}}\atop
{\textstyle\smash{\mathop{\longleftarrow}\limits_{\rho_k^N}}}&\ca O_S&\longrightarrow&O\cr
&&\lmapdown{\psi}&&\lmapdown{\psi}&&\Big\|&\cr
O&\longrightarrow&\ca I_S/\ca I_S^{k+1}&
\longhookrightarrow&\ca O_M/\ca I_S^{k+1}&
{\textstyle\smash{\mathop{\longrightarrow}\limits^{\theta_k}}}\atop
{\textstyle\smash{\mathop{\longleftarrow}\limits_{\rho_k}}}&\ca O_S&\longrightarrow&O\cr
},\cr
\noalign{\bigskip\medskip}
\matrix{O&\longrightarrow&\ca
I_{S,N}^k/\ca I_{S,N}^{k+1}&\longhookrightarrow&\ca I_{S,N}/\ca I_{S,N}^{k+1}&
{\textstyle\smash{\mathop{\longrightarrow}\limits^{\theta_{k,k-1}^N}}}\atop
{\textstyle\smash{\mathop{\longleftarrow}\limits_{\nu_{k-1,k}^N}}}&\ca
I_{S,N}/\ca I_{S,N}^k&\longrightarrow&O\cr
&&\lmapdown{\psi}&&\lmapdown{\psi}&&\lmapdown{\hat\psi}&\cr
O&\longrightarrow&\ca I_S^k/\ca I_S^{k+1}&
\longhookrightarrow&\ca I_S/\ca I_S^{k+1}&
{\textstyle\smash{\mathop{\longrightarrow}\limits^{\theta_{k,k-1}}}}\atop
{\textstyle\smash{\mathop{\longleftarrow}\limits_{\nu_{k-1,k}}}}&\ca I_S/\ca I_S^k&
\longrightarrow&O\cr
}.\cr}
$$

First of all, we define $\rho_k=\psi\circ\rho_k^N$. As in the proof of
Proposition~\usette, we see that this is a ring morphism such that
$\theta_k\circ\rho_k=\id$,
and so $S$ is $k$-split in~$M$.

Now, $\theta_k\circ\psi=\theta_k^N$ implies that $\psi(\hbox{\rm Ker }\theta_k^N)\subseteq
\hbox{\rm Ker }\theta_k$, that is $\psi(\ca I_{S,N}/\ca I_{S,N}^{k+1})\subseteq \ca I_S/\ca
I_S^{k+1}$. Since $\psi$ is a ring isomorphism, it induces a ring isomorphism (still
denoted by~$\psi$) between the nilpotent parts of the two rings,~$\ca I_{S,N}/\ca I_{S,N}^{k+1}$
and~$\ca I_S/\ca I_S^{k+1}$, and thus, by restriction, a ring isomorphism between~$\ca
I_{S,N}^k/\ca I_{S,N}^{k+1}$ and~$\ca I_S^k/\ca I_S^{k+1}$.  Therefore it also induces a quotient
ring isomorphism~$\hat\psi$ between $\ca O_{N_S}/\ca I_{S,N}^k$ and $\ca O_M/\ca I_S^k$ such that
$\theta_{k-1}\circ\hat\psi=\theta^N_{k-1}$, and thus $S_N(k-1)\cong S(k-1)$. Furthermore,
$\hat\psi$ sends~$\ca I_{S,N}/\ca I_{S,N}^k$ into~$\ca I_S/\ca I_S^k$ so that
$\hat\psi\circ\theta_{k,k-1}^N=\theta_{k,k-1}\circ\psi$. This restriction of~$\hat\psi$ is
an isomorphism of~$\ca O_S$-modules: indeed
$$
\eqalign{
\hat\psi(u\cdot[h]_k)&=\hat\psi\bigl(\rho_{k-1}^N(u)[h]_k\bigr)=\hat\psi\bigl(
\theta_{k,k-1}^N(\rho_k^N(u)[h]_{k+1})\bigr)=\bigl[\psi\bigl(\rho_k^N(u)\bigr)\psi[h]_{k+1}
\bigr]_k\cr
&=\bigl[\rho_k(u)\psi[h]_{k+1}\bigr]_k=\rho_{k-1}(u)\hat\psi([h]_k)=u\cdot\hat\psi([h]_k)
\cr}
$$
for all $u\in\ca O_S$ and $h\in\ca I_{S,N}$, where
$\rho_{k-1}^N=\theta_{k,k-1}^N\circ\rho_k^N$ and $\rho_{k-1}=\theta_{k,k-1}\circ\rho_k$.

We then define $\nu_{k-1,k}=\psi\circ\nu_{k-1,k}^N\circ\hat\psi^{-1}$; we claim that
$\nu_{k-1,k}$ is a morphism of~$\ca O_S$-modules such that
$\theta_{k,k-1}\circ\nu_{k-1,k}=\id$. Indeed,
$$
\eqalign{
\nu_{k-1,k}(u\cdot[h]_k)&=\psi\circ\nu_{k-1,k}^N\circ\hat\psi^{-1}(u\cdot[h]_k)=
\psi\circ\nu_{k-1,k}^N\bigl(u\cdot\hat\psi^{-1}([h]_k)\bigr)=\psi\bigl(u\cdot
(\nu_{k-1,k}^N\circ\hat\psi^{-1})([h]_k)\bigr)\cr
&=\psi\bigl(\rho_k^N(u)(\nu_{k-1,k}^N\circ\hat\psi^{-1})([h]_k)\bigr)=
\psi\bigl(\rho_k^N(u)\bigr)(\psi\circ\nu_{k-1,k}^N\circ\hat\psi^{-1})([h]_k)
=\rho_k(u)\nu_{k-1,k}([h]_k)\cr
&=u\cdot\nu_{k-1,k}([h]_k)
\cr}
$$
for all $u\in\ca O_S$ and $h\in\ca I_S$. Finally,
$\theta_{k,k-1}\circ\nu_{k-1,k}=\hat\psi\circ\theta^N_{k,k-1}
\circ\nu_{k-1,k}^N\circ\hat\psi^{-1}=\id$, and hence $S$ is
$(k-1)$-comfortably embedded in~$M$, as claimed.

Conversely, assume that $S$ is $k$-split and $(k-1)$-comfortably embedded.
Since we shall use different maps, let us write the involved commutative diagrams:
$$

\displaylines{\matrix{
O&\longrightarrow&\ca
I_{S,N}/\ca I_{S,N}^{k+1}&\longhookrightarrow&\ca O_{N_S}/\ca I_{S,N}^{k+1}&
{\textstyle\smash{\mathop{\longrightarrow}\limits^{\theta_k^N}}}\atop
{\textstyle\smash{\mathop{\longleftarrow}\limits_{\rho_k^N}}}&\ca O_S&\longrightarrow&O\cr
&&\lmapdown{\psi^k}&&\lmapdown{\psi}&&\Big\|&\cr
O&\longrightarrow&\ca I_S/\ca I_S^{k+1}&
\longhookrightarrow&\ca O_M/\ca I_S^{k+1}&
{\textstyle\smash{\mathop{\longrightarrow}\limits^{\theta_k}}}\atop
{\textstyle\smash{\mathop{\longleftarrow}\limits_{\rho_k}}}&\ca O_S&\longrightarrow&O\cr
},\cr
\noalign{\bigskip\medskip}
\matrix{O&\longrightarrow&\ca
I_{S,N}^l/\ca
I_{S,N}^{l+1}&{\textstyle\smash{\longhookrightarrow}}\atop
{\textstyle\smash{\mathop{\longleftarrow}\limits_{\tilde\nu_{l-1,l}^N}}}&\ca I_{S,N}/\ca
I_{S,N}^{l+1}& {\textstyle\smash{\mathop{\longrightarrow}\limits^{\theta_{l,l-1}^N}}}\atop
{\textstyle\smash{\mathop{\longleftarrow}\limits_{\nu_{l-1,l}^N}}}&\ca I_{S,N}/\ca
I_{S,N}^l&\longrightarrow&O\cr
&&\lmapdown{\chi^l}&&\lmapdown{\psi^l}&&\lmapdown{\psi^{l-1}}&\cr O&\longrightarrow&\ca
I_S^l/\ca I_S^{l+1}&
{\textstyle\smash{\longhookrightarrow}}\atop
{\textstyle\smash{\mathop{\longleftarrow}\limits_{\tilde\nu_{l-1,l}^N}}}&\ca I_S/\ca
I_S^{l+1}&
{\textstyle\smash{\mathop{\longrightarrow}\limits^{\theta_{l,l-1}}}}\atop
{\textstyle\smash{\mathop{\longleftarrow}\limits_{\nu_{l-1,l}}}}&\ca I_S/\ca I_S^l&
\longrightarrow&O\cr
}.\cr}
$$

In the proof of Proposition~\usette\ we defined an $\ca O_S$-module isomorphism
$\chi\colon\ca I_{S,N}/\ca I_{S,N}^2\to\ca I_S/\ca I_S^2$; since $\ca I_S^l/\ca
I_S^{l+1}=\hbox{\rm Sym}^l(\ca I_S/\ca I_S^2)$, and likewise for $\ca I_{S,N}^l/\ca
I_{S,N}^{l+1}$, we get for all $l\ge 1$ an $\ca O_S$-module isomorphism $\chi^l\colon\ca
I_{S,N}^l/\ca I_{S,N}^{l+1}\to\ca I_S^l/\ca I_S^{l+1}$. We claim that for all $1\le l\le k$
we can define a ring and $\ca O_S$-module isomorphism $\psi^l\colon\ca I_{S,N}/\ca
I_{S,N}^{l+1}\to\ca I_S/\ca I_S^{l+1}$ so that the above diagram commutes.

We argue by induction on~$l$. For $l=1$, it suffices to take $\psi^1=\chi$. Assume now
that we have defined~$\psi^{l-1}$, and let
$$
\psi^l=\nu_{l-1,l}\circ\psi^{l-1}\circ\theta_{l,l-1}^N+\chi^l\circ\tilde\nu_{l-1,l}^N,
$$
where $\tilde\nu_{l-1,l}^N$ is the left splitting morphism associated to~$\nu_{l-1,l}^N$.
It is easy to check that~$\psi^l$ is invertible, with inverse given by
$$
(\psi^l)^{-1}=\nu^N_{l-1,l}\circ(\psi^{l-1})^{-1}\circ\theta_{l,l-1}+(\chi^l)^{-1}
\circ\tilde\nu_{l-1,l},
$$
where $\tilde\nu_{l-1,l}$ is the left splitting morphism associated to~$\nu_{l-1,l}$.
Since all the maps involved are $\ca O_S$-module morphisms, $\psi^l$ is a $\ca O_S$-module
morphism; we are left to show that $\psi^l$ is a ring morphism. Using the definition, it is
easy to see that $\psi^l$ is a ring morphism if and only if
$$
\tilde\nu_{l-1,l}\bigl((\nu_{l-1,l}\circ\psi^{l-1})(u)(\nu_{l-1,l}\circ\psi^{l-1})
(v)\bigr)=\chi^l\circ\tilde\nu_{l-1,l}^N
\bigl(\nu^N_{l-1,l}(u)\nu_{l-1,l}^N(v)\bigr),
\neweq\eqprodo
$$
for all $u$,~$v\in\ca I_{S,N}/\ca I_{S,N}^l$.

To prove \eqprodo, we work in local coordinates. Let $(U,z)$ be a chart in a
$(k-1)$-comfortable atlas, so that
we have $\nu_{i-1,i}[z^{r_1}\cdots z^{r_a}]_i=[z^{r_1}\cdots z^{r_a}]_{i+1}$ for all
$i=1,\ldots,l$, $a=1,\ldots,l-1$ and $r_1,\ldots,r_a=1,\ldots,m$. The chart $(U,z)$ induces
a chart $(\tilde U,\tilde z)$ in a $(k-1)$-comfortable atlas of~$N_S$, with
$\tilde z=(v,z'')$, where $v=(v^1,\ldots,v^m)$ are the fiber coordinates. Then we have
$\chi[v^r]_2=[z^r]_2$, and thus
$$
\chi^i[v^{r_1}\cdots v^{r_i}]_{i+1}=[z^{r_1}\cdots  z^{r_i}]_{i+1}
$$
for all $i=1,\ldots,k$. From this and the fact that
$\tilde\nu_{i-1,i}^N([v^r]_{i+1})=O$ for all~$i=1,\ldots,k$ and~$r=1,\ldots,m$, it follows
easily that $\psi^i([v^r]_{i+1})=[z^r]_{i+1}$ for all $i=1,\ldots,l-1$ and~$r=1,\ldots,m$.

Now take $u$,~$v\in\ca I_{S,N}/\ca I_{S,N}^l$; we can write
$$
u=\sum_{a=1}^{l-1} \alpha_{r_1\ldots r_a}\cdot[v^{r_1}\cdots
v^{r_a}]_l\qquad\hbox{and}\qquad
u=\sum_{b=1}^{l-1} \beta_{s_1\ldots s_b}\cdot[v^{s_1}\cdots v^{s_b}]_l
$$
for suitable $\alpha_{r_1\ldots r_a}$,~$\beta_{s_1\ldots s_b}\in\ca O_S$. Then
$$
\chi^l\circ\tilde\nu_{l-1,l}^N
\bigl(\nu^N_{l-1,l}(u)\nu_{l-1,l}^N(v)\bigr)=\sum_{a+b=l}\alpha_{r_1\ldots r_a}
\beta_{s_1\ldots s_b}\cdot[z^{r_1}\cdots z^{r_a}z^{s_1}\cdots z^{s_b}]_{l+1}.
$$
Using the fact that $\psi^{l-1}$ is a ring morphism and an $\ca O_S$-morphism we also find
$$
(\nu_{l-1,l}\circ\psi^{l-1})(u)=\sum_{a=1}^{l-1} \alpha_{r_1\ldots r_a}\cdot[z^{r_1}\cdots
z^{r_a}]_{l+1}\qquad\hbox{and}\qquad (\nu_{l-1,l}\circ\psi^{l-1})(v)=
\sum_{b=1}^{l-1} \beta_{s_1\ldots s_b}\cdot[z^{s_1}\cdots z^{s_b}]_{l+1},
$$
and hence
$$
\tilde\nu_{l-1,l}\bigl((\nu_{l-1,l}\circ\psi^{l-1})(u)(\nu_{l-1,l}\circ\psi^{l-1})
(v)\bigr)=\sum_{a+b=l}\alpha_{r_1\ldots r_a}
\beta_{s_1\ldots s_b}\cdot[z^{r_1}\cdots z^{r_a}z^{s_1}\cdots z^{s_b}]_{l+1},
$$
as claimed.

So in particular we have proved that $\psi^k\colon\ca I_{N,S}/\ca I_{N,S}^{k+1}\to\ca I_S/
\ca I_S^{k+1}$ is a ring and $\ca O_S$-module isomorphism. Let us then define $\psi\colon\ca
O_{N_S}/\ca I_{S,N}^{k+1}\to\ca O_M/\ca I_S^{k+1}$ by
$$
\psi=\rho_k\circ\theta_k^N+\psi^k\circ \tilde\rho_k^N,
$$
where as usual $\tilde\rho_k^N$ is the derivation associated to~$\rho_k^N$. It is easy to check
that
$\theta_k\circ\psi=\theta_k^N$, and that $\psi$ is invertible; we are left to show that
it is a ring morphism.

If $u$,~$v\in\ca O_{N_S}/\ca I_{S,N}^{k+1}$ we can write
$u=\rho_k^N(u_o)+\tilde\rho_k^N(u)$,
with $u_o=\theta_k^N(u)$; so
$\psi(u)=\rho_k(u_o)+\psi^k\bigl(\tilde\rho_k^N(u)\bigr)$; and analogously for~$v$.
Therefore
$$
\eqalign{
uv&=\rho_k^N(u_o)\rho_k^N(v_o)+\bigl[\rho_k^N(u_o)\tilde\rho_k^N(v)+\rho_k^N(v_o)\tilde\rho_
k^N(u)
+\tilde\rho_k^N(u)\tilde\rho_k^N(v)\bigr]\cr
&=\rho_k^N(u_o)\rho_k^N(v_o)+\bigl[u_o\cdot \tilde\rho_k^N(v)+v_o\cdot \tilde\rho_k^N(u)
+\tilde\rho_k^N(u)\tilde\rho_k^N(v)\bigr],
\cr}
$$
and thus
$$
\eqalign{
\psi(uv)&=\rho_k(u_ov_o)+\psi^k\bigl(u_o\cdot \tilde\rho_k^N(v)+v_o\cdot \tilde\rho_k^N(u)
+\tilde\rho_k^N(u)\tilde\rho_k^N(v)\bigr)\cr
&=\rho_k(u_o)\rho_k(v_o)+u_o\cdot\psi^k\bigl(\tilde\rho_k^N(v)\bigr)+v_o\cdot\psi^k\bigl(
\tilde\rho_k^N(u)\bigr)+\psi^k\bigl(\tilde\rho_k^N(u)\bigr)\psi^k\bigl(\tilde\rho_k^N(v)
\bigr)\cr
&=\psi(u)\psi(v).
\cr}
$$
\qedn

Thus a submanifold $S$ is $k$-linearizable if and only if there is an atlas
$\{(U_\alpha,z_\alpha)\}$ adapted to~$S$ whose changes of coordinates are of the form
$$
\cases{z_\beta^r=(a_{\beta\alpha})^r_s(z''_\alpha) z_\alpha^s+R_{k+1}&for $r=1,\ldots,m$,\cr
\noalign{\smallskip}
z_\beta^p=\phi^p_{\alpha\beta}(z''_\alpha)+R_{k+1}&for $p=m+1,\ldots,n$.\cr}
$$

\Rem Camacho, Movasati and Sad in [CMS] defined a $k$-linearizable curve as a complex
curve~$S$ in a complex manifold~$M$ for which there exists an atlas adapted to~$S$
whose changes of coordinates are of the form
$$
\cases{z_\beta^r=(a_{\beta\alpha})^r_s(z''_\alpha) z_\alpha^s+R_{k+1}&for $r=1$,\cr
\noalign{\smallskip}
z_\beta^p=\phi^p_{\alpha\beta}(z''_\alpha)&for $p=2,\ldots,n$;\cr}
$$
they dropped the remainder term in the $z''_\alpha$ variables only because they were interested
in curves with a neighbourhood fibered by $(n-1)$-dimensional disks. As a consequence, our notion
of 2-linearizable curves (or submanifolds) is strictly weaker than the notion of
2-linearizable curves used in [CMS].

Recalling Proposition~\ksplit\ and Corollary~\obsce\ we see that the obstructions for passing
from $(k-1)$-linearizable to $k$-linearizable live in 
$H^1\bigl(S,\ca T_S\otimes\hbox{Sym}^k(\ca N_S^*)\bigr)$ and, for $k\ge 2$, in~$H^1\bigl(S,\ca
N_S\otimes\hbox{Sym}^k(\ca N_S^*)\bigr)$. Using again Grauert's vanishing theorem [G,
Hilfssatz~1, p. 344] we get

\newthm Corollary \tutlin: Let $S$ be an $m$-codimensional compact complex submanifold of an
$n$-dimensional manifold~$M$, and assume that $N_S$ is negative in the sense of Grauert. Then there
exists a~$k_0\ge 1$ such that if~$S$ is $k_0$-linearizable then it is $k$-linearizable for all
$k\ge k_0$.

Again, we can get similar results also assuming suitable positivity conditions on~$N_S$; see~[Gr,
CG, K1, K2, St]. Furthermore, in the next section we shall be able to compute the number $k_0$ for
curves in a complex surface. 

\Rem
When $S$ is a hypersurface in~$M$, and thus $N_S$ is a line bundle, we actually have found that
the obstructions to $k$-linearizability live in $H^1\bigl(S,\ca
T_S\otimes(\ca N_S^*)^{\otimes k}\bigr)$ and in~$H^1\bigl(S,\ca (\ca N_S^*)^{\otimes k-1}\bigr)$,
in accord with Grauert's theory (see, again, [G] and [CM]).

\Rem
It is important to remark that the isomorphisms between $S(k)$ and $S_N(k)$ obtained in the
previous corollary are compatible in the sense that if $k'>k$ then the restriction to $S(k)$ of
the isomorphism between $S(k')$ and $S_N(k')$ induces the given isomorphism between $S(k)$
and~$S_N(k)$. In some sense, we have obtained an isomorphism between the formal neighbourhood
of~$S$ in~$M$ and the formal neighbourhood of~$S$ in~$N_S$. Grauert [G] and others have given
conditions ensuring that such a formal isomorphism extends to a biholomorphism between an actual
neighbourhood of $S$ in~$M$ to an actual neighbourhood of $S$ in~$N_S$. In particular, applying
Grauert's formal principle (see [CM, Theorem~4.3]) we recover Grauert's result:

\newthm Corollary \Grauert: (Grauert [G]) Let $S$ be a compact complex hypersurface of an
$n$-dimensional manifold~$M$. Assume that $N_S$ is negative in the sense of Grauert, and that $S$
is exceptional in~$M$ (that is, it can be blown down to a point). Then there exists a $k_0\ge 1$
such that if $S$ is $k_0$-linearizable then a neighbourhood of $S$ in~$M$ is biholomorphic to a
neighbourhood of~$S$ in~$N_S$.

\Rem
A compact Riemann surface $S$ in a complex surface $M$ is exceptional if and only if its
self-intersection $S\cdot S$ is negative, if and only if its normal bundle is negative in the sense
of Grauert; see [CM].

\Rem
The formal principle holds in several other instances too (but not always). For instance, using
[Gr, CG, K1, K2, St] we can get a statement analogous to Corollary~\Grauert\ assuming suitable
positivity conditions on~$N_S$ (and arbitrary codimension).

%
%
%
%
%
%
\Rem Let  $S$ be a submanifold of a complex manifold $M$ and denote by $g\colon S\to M$ the
inclusion. Then one of these cases hold:
\smallskip

\item{a)} $S$ is $k$-linearizable in $M$ for all $k\geq 1$.

\item{b)} $S$ is $k$-split in $M$ for all $k\ge 1$, but there exist $m_g\geq 1$ and a
$m_g$-splitting $\rho$ such that $S$ is not $m_g$-comfortably embedded in~$M$ but it is
$(m_g-1)$-comfortably embedded. In this case, we can
associate to the embedding $S\to M$ a non zero
cohomology class  $\gt h^{\rho}\in H^1\bigl(S,\ca N_S\otimes
\ca I_S^{m_g+1}/\ca I_S^{m_g+2}\bigr)$.

\item{c)} there exists an integer $k_g\geq 1$ such that $S$ is not  $k$-split in $M$, but it is
$(k-1)$-splitting. In this case, given any fixed $(k_g-1)$-splitting $\rho$, we can associate to
the embedding $S\to M$ a non zero cohomology
class $\gt g_{k_g} \in H^1\bigl(S,
\ca T_S\otimes\ca I_S^{k_g}/\ca I_S^{k_g+1}\bigr)$. Furthermore, we can choose $1\le m_g\le k_g-1$
so that $S$ is $(m_g-1)$-comfortably embedded (with respect to the lifting induced by $\rho$)
in~$M$ but not $m_g$-comfortably embedded in~$M$, and hence we get a non zero
cohomology class  $\gt h^{\rho}\in H^1\bigl(S,\ca N_S\otimes
\ca I_S^{m_g+1}/\ca I_S^{m_g+2}\bigr)$.
\smallskip
\noindent It is clear by the construction that if two different embeddings of the same
submanifold~$S$ have biholomorphic neighbourhoods then the integers and the cohomology classes
constructed above must be the same in both cases. It would be interesting to know other
invariants. For instance, a consequence of Corollary~\Grauert\ and Remark~4.4 is that two
infinitely linearizable (in the sense of Remark~4.3) embeddings of a compact Riemann surface with
negative self-intersection in a complex surface always have biholomorphic neighbourhoods.
However, as far as we know, even for curves with negative self-intersection in case (b) other
invariants beside~$m_g$ and $\gt h^\rho$ are not yet known.

%
%
\smallsect 5. Embeddings of a smooth curve

In this section we shall use Serre duality to describe sufficient conditions for a compact curve in
a complex surface to be $k$-split, $k$-comfortably embedded and/or $k$-linearizable.

%
%
%

Let $S$ be a non-singular, compact, irreducible curve of genus~$g$ on a surface~$M$. In
particular, $N_S$ is a line bundle; therefore $\hbox{Sym}^k(\ca N_S^*)\cong (\ca N_S^*)^{\otimes
k}$ for all $k\ge 1$, and the obstruction for passing from
$(k-1)$-split to
$k$-split lives in
$H^1\bigl(S,\ca T_S\otimes(\ca N_S^*)^{\otimes k}\bigr)$.
The
Serre duality for Riemann surfaces implies that
$$
H^1\bigl(S,\ca T_S\otimes(\ca N_S^*)^{\otimes k}\bigr)\cong
H^0\bigl(S,\Omega_S\otimes\Omega_S\otimes\ca N_S^{\otimes k}\bigr).
$$
Now,
$$
\deg(\Omega_S\otimes\Omega_S\otimes {\ca N}_S^{\otimes k})=4g-4+k(S\cdot S);
$$
therefore 
$$
k(S\cdot S)<4-4g\qquad\Longrightarrow\qquad H^1\bigl(S,\ca T_S\otimes(\ca N_S^*)^{\otimes k}\bigr)
=(O).
\neweq\eqksplitcurv
$$ 
It follows in particular that if $g\ge 1$ and $S\cdot S<4-4g$, or $g=0$ and $S\cdot S\le 0$, then
$S$ is $k$-splitting in~$M$ for every~$k\ge 1$.

The obstruction for a split curve to be 1-comfortably embedded is in $H^1\bigl(S,\ca
N_S\otimes(\ca N_S^*)^{\otimes 2}\bigr)\cong H^1(S,\ca N_S^*)$. Serre duality yields
$$
H^1(S,\ca N_S^*)\cong H^0(S,\Omega_S\otimes\ca N_S);
$$
so, since $\deg(\Omega_S\otimes {\ca N}_S)=2g-2+S\cdot S$, we get
$$
S\cdot S<2-2g \qquad\Longrightarrow\qquad H^1\bigl(S,\ca
N_S\otimes(\ca N_S^*)^{\otimes 2}\bigr)=(O).
\neweq\eqkcecurv
$$ 
In particular, if
$g\ge 1$ and
$S\cdot S<4-4g$ or $g=0$ and $S\cdot S<2$ then $S$ is (splitting and) $1$-comfortably embedded. 

More generally, assume that $S$ is $k$-split and $(k-1)$-comfortably embedded. The obstruction for
$S$ to be
$k$-comfortably embedded lives in $H^1\bigl(S,\ca
N_S\otimes(\ca N_S^*)^{\otimes k+1}\bigr)\cong H^1\bigl(S,(\ca N_S^*)^{\otimes k})$. Then using
Serre duality as before we find
$$
k(S\cdot S)< 2-2g \qquad\Longrightarrow\qquad H^1\bigl(S,\ca
N_S\otimes(\ca N_S^*)^{\otimes k+1}\bigr)=(O).
\neweq\eqkcecurvk
$$
%
In particular, if $g\ge 1$ and $S\cdot
S<2-2g$ then $k$-splitting implies $k$-comfortably embedded, while if $g=0$ and
$S\cdot S\le 0$ then $S$ is $k$-comfortably embedded for all $k\ge 1$.

We can summarize the content of our computations in the following

\newthm Proposition \cusu: Let $S$ be a non-singular, compact, irreducible curve of
genus~$g$ in a surface~$M$. Then:
\smallskip
\itm{(i)}if $g\ge 1$ and $S\cdot S<4-4g$ then $S$ is $k$-split into~$M$ for all $k\ge 1$;
\item{\rm (ii)}if $g\ge 1$ and $S\cdot S<2-2g$ then $S$ $k$-split implies $S$ $k$-comfortably
embedded into~$M$ for any~$k\ge 1$; in particular, if $g\ge 1$ and $S\cdot S<4-4g$ then $S$
is $k$-linearizable for all $k\ge 1$;
\itm{(iii)}if $g=0$ and $S\cdot S\le 0$ then $S$ is $k$-linearizable for all $k\ge 1$;
\itm{(iv)} if $g=0$ and $S\cdot S\le 1$ then $S$ is $3$-split and $1$-comfortably embedded
into~$M$;
\itm{(v)} if $g=0$ and $S\cdot S\le 3$ then $S$ splits into~$M$.

\Rem 
Proposition~\cusu.(ii) has been proved in a slightly different way in [CMS],
where $S$ was assumed to be fibered imbedded into~$M$ (and thus, in particular, $k$-split for
all $k\ge 1$).

Another way of looking at \eqksplitcurv\ and \eqkcecurvk\ yields the following

\newthm Proposition \cusb: Let $S$ be a non-singular, compact, irreducible curve of genus~$g\ge 1$
in a surface~$M$. Then:
\smallskip
\itm{(i)}if $S\cdot S<0$ and $S$ is $k_0$-splitting for some $k_0> (4g-4)/|S\cdot S|$ then $S$
is $k$-splitting for all $k\ge k_0$;
\itm{(ii)}if $S\cdot S<0$ and $S$ is $k_0$-comfortably embedded for some $k_0> (2g-2)/|S\cdot
S|$ then $k$-splitting implies\break\indent $k$-comfortably embedded for any $k\ge k_0$;
\itm{(iii)}if $S\cdot S<0$ and $S$ is $k_0$-linearizable for some $k_0> (4g-4)/|S\cdot S|$ then $S$
is $k$-linearizable for all $k\ge k_0$.

Recalling Remarks~4.3 and~4.4, we can apply Grauert's formal principle ([CM, Theorem~4.3]) to
recover, among other things, results due to Laufer and Camacho-Movasati-Sad:

\newthm Corollary \Laufer: Let $S$ be a non-singular, compact, irreducible curve
of genus~$g$ in a surface~$M$ with negative self-intersection~$S\cdot S<0$. If
\smallskip
\itm{(a)} $g=0$, or
\item{\rm(b)} $g\ge 1$, $S$ is $k_0$-split and $k_1$-comfortably embedded for some $k_0>
(4g-4)/|S\cdot S|$ and $k_1>(2g-2)/|S\cdot S|$, or
\itm{(c)} {\rm (Laufer [L, Chapter V\negthinspace I])} $g\ge 1$ and $S\cdot S<4-4g$, or
\itm{(d)} {\rm [CMS]} $g\ge 1$, $S\cdot S<2-2g$ and $S$ is $k_0$-split for some $k_0>
(4g-4)/|S\cdot S|$,
\smallskip
\noindent then a neighbourhood
of $S$ in~$M$ is biholomorphic to a neighbourhood of the zero section of~$N_S$.

\smallsect 6. Another characterization of split and comfortably embedded submanifolds

In [ABT2] we used the 1-comfortably embedded condition to build partial
holomorphic connections on the normal bundle, and we wondered why this condition appeared to be
the right one for such constructions. In this section we give an answer of sort to this question,
showing that a submanifold is 1-comfortably embedded  if and
only if it exists an infinitesimal holomorphic connection on the normal bundle.

Let us begin with a definition.
 

\Def Let $S$ be a complex subvariety of a complex manifold~$M$. The
{\sl sheaf of holomorphic differentials} on~$S(1)$ is given by
$$
\Omega_{S(1)}=\Omega_M/(\ca I_S^2\Omega_M+d\ca I_S^2);
$$
its dual $\ca T_{S(1)}= {\rm Hom}_{{\ca
O}_{S(1)}} (\Omega_{S(1)}, \ca O_{S(1)})$, where $\ca O_{S(1)}=\ca O_M/\ca I_S^2$ as usual, is the
{\sl holomorphic tangent sheaf} of~$S(1)$. The map $d_{(1)}\colon\ca O_{S(1)}\to\Omega_{S(1)}$
given by $d_{(1)}([f]_2)=\pi_1(df)$, where
$\pi_1\colon\Omega_M\to\Omega_{S(1)}$ is the natural projection, is the {\sl canonical
differential.} We refer to [L, Chapter V\negthinspace I]
for properties of differentials on an analytic space with nilpotents.

%
%

Theorem~\ucinque\ yields a characterization of splitting manifolds in terms of~$\Omega_{S(1)}$:

\newthm Proposition \carinfsplit: Let $S$ be a submanifold of a complex manifold~$M$. 
Then $S$ splits into~$M$ if and only if there exists a surjective 
$\ca O_{S(1)}$-morphism 
$$
X_{(1)}\colon\Omega_{S(1)}\to \ca I_S/\ca I_S^2
$$ 
such that $X_{(1)}\circ d_{(1)}\circ i_1=\id$, where $i_1\colon\ca I_S/\ca I_S^2\to\ca O_M/\ca
I_S^2$  is the natural inclusion. Furthermore, 
if $\rho\colon\ca O_S\to\ca O_{S(1)}$ is a first order lifting,
and
$\tilde\rho\colon\ca O_{S(1)}\to\ca I_S/\ca I_S^2$ is the associated left splitting morphism,
then $\tilde\rho=X_{(1)}\circ d_{(1)}$.


\pf By Theorem~\ucinque, $S$ splits in $M$ if and only if there
exists a $\theta_1$-derivation $\tilde\rho\colon \ca O_{S(1)}\to \ca I_S/\ca I^2_S$ such
that~$\tilde\rho\circ i_1 =\id$. 

Assuming $\tilde\rho$ given, the universal property of differentials
yields a $\ca O_{S(1)}$-morphism
$X_{(1)}\colon\Omega_{S(1)}\to\ca I_S/\ca I_S^2$ such that
 $X_{(1)}\circ
d_{(1)} = \tilde\rho$; in particular, $X_{(1)}\circ
d_{(1)}\circ i_1=\id$. 

Conversely, given $X_{(1)}$ then $\tilde\rho = X_{(1)}\circ d_{(1)}$
is a  $\theta_1$-derivation such that $\tilde\rho\circ i_1=\id$, and thus $S$ splits into~$M$.
\qedn

To give the announced characterization of 1-comfortably embedded submanifolds we need a last
definition and a last proposition.

\Def Let $S$ be a submanifold of a complex manifold~$M$. Assume that $S$ splits in $M$ and
denote  by
 $X_{(1)}\colon\Omega_{S(1)}\to\ca I_S/\ca I_S^2$ the $\ca O_{S(1)}$-morphism associated
to the choice of a first order lifting by  Proposition~\carinfsplit. An {\sl infinitesimal normal
connection along~$X_{(1)}$} on a  
$\ca O_{S(1)}$-module~$\ca E$ on~$S$ is a $\C$-linear map
$\tilde X_{(1)}\colon\ca E\to \ca I_S/\ca I_S^2\otimes_{\ca O_{S(1)}}\ca E$
satisfying the Leibniz rule
$$
\tilde X_{(1)}(gs)=X_{(1)}\bigl(d_{(1)}(g)\bigr)\otimes s+g\tilde X_{(1)}(s)
$$
for all local sections $g$ of $\ca O_{S(1)}$ and $s$ of $\ca
E$.

\Rem Any locally free $\ca O_S$-module $\ca E$ can be considered as a locally free $\ca
O_{S(1)}$-module endowing it with the structure obtained by restriction of the scalars via the
first order lifting~$\rho$. However, for the application we have in mind we shall need a locally
free $\ca O_{S(1)}$-module which is not obtained in this way.

\Rem
In this section, indeces like $a$,~$b$,~$c$,~$d$ will run from~1 to~$\hbox{rk}(\ca E)$.

\newthm Proposition \carnormjet: Let $S$ be a submanifold of a complex manifold~$M$. 
Assume that $S$ splits in
$M$, with first order lifting $\rho\colon {\ca O}_S\to {\ca O}_{S(1)}$ and associated
$\ca O_{S(1)}$-morphism $X_{(1)}\colon\Omega_{S(1)}\to {\ca N}^*_S$.  
Let $\ca E$ be a locally free $\ca O_{S(1)}$-module on~$S$. Then
the obstruction to the existence of an infinitesimal
normal connection on $\ca E$ along~$X_{(1)}$ is the class 
$\delta_\rho({\cal E})\in H^1(S,\ca I_S/\ca I_S^2\otimes\End ({\cal E}))$ represented,
in an atlas $\gt U=\{(U_\alpha,z_\alpha)\}$ adapted to $\rho$ and trivializing~$\ca E$, by the
$1$-cocycle
$$
[(\Phi_{\beta\alpha})^{c}_{a}]_2
\left[{\de(\Phi_{\alpha\beta})^{d}_{c}\over \de z^r_\alpha}
z^r_\alpha\right]_2
\otimes e^{*a}_\beta\otimes e_{d,\beta},$$
where $ e_{b,\alpha}$ (for $b=1,\ldots, {\rm rk}\,\ca E$) is a local frame
for $\ca E$ over~$U_\alpha\cap S$, $e^{* b}_\alpha$ is the dual frame, and 
$[(\Phi_{\alpha\beta})^b_c]_2\in\ca O_{S(1)}$ are the transition functions of~$\ca E$.

\pf 
Let 
$\tilde{X}_{(1)}\colon{\cal E}\to \ca I_S/\ca I_S^2 \otimes_{{\ca O}_S}
{\cal E}$
 be an
infinitesimal normal connection along~$X_{(1)}$, and define an element
$\eta^b_{c,\alpha}\in H^0(U_\alpha\cap S,\ca I_S/\ca I_S^2)$ by the formula
$$
\tilde X_{(1)}(e_{c,\alpha})=\eta^b_{c,\alpha}\otimes e_{b,\alpha}.
$$
Now, if $U_\alpha\cap U_\beta\cap S\ne\void$ we have
$
e_{b,\alpha}=[(\Phi_{\alpha\beta})^d_b]_2   e_{d,\beta};
$
so
$$
\eqalign{\eta^b_{c,\alpha}\otimes [(\Phi_{\alpha\beta})^d_b]_2\ e_{d,\beta}&=\tilde
X_{(1)}(e_{c,\alpha})=\tilde X_{(1)}\left([(\Phi_{\alpha\beta})^d_c]_2\ e_{d,\beta}\right)\cr
&=X_{(1)}\left(d_{(1)}[(\Phi_{\alpha\beta})^d_c]_2\right)\otimes
e_{d,\beta}+[(\Phi_{\alpha\beta})^b_c]_2 \cdot\eta^d_{b,\beta}
\otimes e_{d,\beta}.
\cr}
$$
But 
$$
X_{(1)}\left(d_{(1)}([(\Phi_{\alpha\beta})^d_c]_2)\right)
=\tilde\rho\left([(\Phi_{\alpha\beta})^d_c]_2\right)=
\left[{\de (\Phi_{\alpha\beta})^d_c\over\de z_\alpha^r} z_\alpha^r\right]_2,
$$
by Remark~1.8, and hence
$$
[(\Phi_{\alpha\beta})^d_b]_2\cdot \eta^b_{c,\alpha}=
\left[{\de (\Phi_{\alpha\beta})^d_c\over\de z_\alpha^r} z_\alpha^r\right]_2+
[(\Phi_{\alpha\beta})^b_c]_2 \cdot\eta^d_{b,\beta}.
$$
If we define the 0-cocycle $\gt k=\{\gt k_\alpha\}\in H^0(\gt U_S,\ca I_S/\ca I_S^2\otimes
\ca E^*\otimes\ca E)$ by setting
$$
\gt k_\alpha=\eta^b_{c,\alpha}\otimes e^{*c}_\alpha\otimes e_{b,\alpha},
$$
we get
$$
\eqalign{
\gt k_\alpha-\gt k_\beta&= \eta^b_{c,\alpha}\otimes e^{*c}_\alpha\otimes e_{b,\alpha}
- \eta^b_{c,\beta}\otimes e^{*c}_\beta\otimes e_{b,\beta}\cr
&=
[(\Phi_{\beta\alpha})^b_d]_2\left(\left[{\de (\Phi_{\alpha\beta})^d_c\over\de z_\alpha^r}
z_\alpha^r\right]_2+ [(\Phi_{\alpha\beta})^a_c]_2 \cdot\eta^d_{a,\beta}\right)
\otimes [(\Phi_{\beta\alpha})^c_d]_2 e^{*d}_\beta\otimes [(\Phi_{\alpha\beta})^d_b]_2 \ 
e_{d,\beta} - \eta^b_{c,\beta}\otimes e^{*c}_\beta\otimes e_{b,\beta}
\cr
&=
\left[{\de (\Phi_{\alpha\beta})^d_c\over\de z_\alpha^r}
z_\alpha^r\right]_2
\otimes [(\Phi_{\beta\alpha})^c_a]_2 e^{*a}_\beta\otimes \ 
e_{d,\beta},
\cr
}
$$
and thus $\delta_\rho(\ca E)=O$.

Conversely, assume that $[(\Phi_{\beta\alpha})^{c}_{a}]_2
\left[{\de(\Phi_{\beta
\alpha})^{d}_{c}\over \de z^r_\alpha}\right]_2
\left[z^r_\alpha\right]_2
\otimes e^{*a}_\beta\otimes e_{d,\beta}=\gt k_\alpha-\gt k_\beta$
with $\gt k_\alpha\in H^0(\gt U_S,\ca I_S/\ca I_S^2\otimes
\ca E^*\otimes\ca E)$. Writing $
\gt k_\alpha=\eta^b_{c,\alpha}\otimes e^{\ast,c}_\alpha\otimes e_{b,\alpha}$, it is easy to
check that setting
$$
\tilde X_{(1)}(e_{c,\alpha})=\eta^b_{c,\alpha}\otimes e_{b,\alpha}
$$
we define an infinitesimal normal connection on $\ca E$.
\qedn

If $S$ splits into~$M$, and $\gt U=\{(U_\alpha,z_\alpha)\}$ is a splitting atlas, then it is easy
to check that the position
$$
(\Phi_{\alpha\beta})^r_s=\left[{\de z_\beta^r\over\de z_\alpha^s}\right]_2
$$
defines a 1-cocycle with coefficients in $GL(m,\ca O_{S(1)})$, and hence a locally free $\ca
O_{S(1)}$-module on~$S$ that, with a slight abuse of notations, we shall denote by~$\ca N_S$. One
of the reasons justifying this notation is that the 1-cocycle of the locally free
$\ca O_{S(1)}$-module $\ca I_S/\ca I_S^2$ is the inverse transposed of the 1-cocycle of~$\ca N_S$,
and thus~$\ca N_S^*\cong \ca I_S/\ca I_S^2$ as $\ca O_{S(1)}$-modules too. Notice, however, that
\defgrho\ implies that this~$\ca N_S$ is not the locally free $\ca O_{S(1)}$-module obtained by
restriction of the scalars via~$\rho$ starting from the usual normal sheaf on~$S$ (as described in
Remark~6.1) {\it unless} $S$ is 1-comfortably embedded in~$M$. 

We finally have the promised characterization of 1-comfortably embedded submanifolds:

\newthm Proposition \infcomftemb: Let $S$ be a submanifold of a complex manifold~$M$. 
Assume that $S$ splits in
$M$, with first order lifting $\rho\colon {\ca O}_S\to {\ca O}_{S(1)}$ and associated
$\ca O_{S(1)}$-morphism $X_{(1)}\colon\Omega_{S(1)}\to \ca I_S/\ca I_S^2$. Then
 $S$ is $1$-comfortably embedded into~$M$ if and only if 
 there exists an infinitesimal normal connection on~$\ca N_S$.

\pf Let $\gt U=\{(U_\alpha,z_\alpha)\}$ be an atlas adapted to~$\rho$, and denote
by~$\{\de_{r,\alpha}\}$ and $\{[z_\alpha^r]_2\}$ the induced local frames on~$\ca N_S$ and~$\ca
I_S/\ca I_S^2$ as locally free $\ca O_{S(1)}$-modules. Proposition~\carnormjet\ says that there
exists an infinitesimal holomorphic connection on~$\ca N_S$ along~$X_{(1)}$ if and only if the
1-cocycle $\delta_\rho(\ca N_S)$ in $H^1\bigl(S,(\ca I_S/\ca I_S^2)^{\otimes 2}\otimes\ca
N_S\bigr)$ given by
$$
\left[{\de z_\alpha^r\over\de z_\beta^s}\right]_2\left[{\de^2 z_\beta^t\over
\de z_\alpha^u\de z_\alpha^r}\right]_2 [z_\alpha^u]_2\otimes [z_\beta^s]_2\otimes\de_{t,\beta}=
\left[{\de^2 z_\beta^t\over
\de z_\alpha^u\de z_\alpha^v}\right]_1[z_\alpha^u]_2\otimes[z_\alpha^v]_2\otimes\de_{t,\beta}
$$
vanishes.

Now, $\delta_\rho(\ca N_S)$ clearly belongs to $H^1\bigl(S,\hbox{Sym}^2(\ca I_S/\ca
I_S^2)\otimes\ca N_S\bigr)$. Since $\hbox{Sym}^2(\ca I_S/\ca I_S^2)$ is a direct summand of~$(\ca
I_S/\ca I_S^2)^{\otimes 2}$, a 1-cocycle in $H^1\bigl(S,\hbox{Sym}^2(\ca I_S/\ca I_S^2)\otimes\ca
N_S\bigr)$ vanishes in $H^1\bigl(S,(\ca I_S/\ca I_S^2)^{\otimes 2}\otimes\ca
N_S\bigr)$ if and only if it vanishes in $H^1\bigl(S,\hbox{Sym}^2(\ca I_S/\ca I_S^2)\otimes\ca
N_S\bigr)$. Since $\hbox{Sym}^2(\ca I_S/\ca I_S^2)\cong \ca I_S^2/\ca I_S^3$, the assertion
follows from Corollary~\obsce.\qedn

\Rem Assume $S$ splits into $M$ and let $X_{(1)}\colon\Omega_{S(1)}\to\ca N_S^*$ be the
corresponding $\ca O_{S(1)}$-morphism. Then one can adapt the notion (see [At]) of first jet
bundle and associate to any $\ca O_{S(1)}$-module~$\ca E$ an $\ca O_{S(1)}$-module $J^1_{\ca
N_S}\ca E$, the {\sl sheaf of normal first jets} of~$\ca E$, and an exact sequence of $\ca
O_{S(1)}$-modules
$$
O\longrightarrow \ca N_S^*\otimes\ca E\longrightarrow J^1_{\ca N_S}\ca E\longrightarrow\ca E
\longrightarrow O
\neweq\eqancheques
$$
in such a way that if $\ca E$ is locally free than $J^1_{\ca N_S}\ca E$ is locally free too, and
the class $\delta_\rho(\ca E)$ introduced in Proposition~\carnormjet\ is exactly the class
associated to the extension~\eqancheques. In particular, Proposition~\infcomftemb\ implies that
the sequence~\eqancheques\ splits if and only if $S$ is 1-comfortably embedded.

%

\bigbreak
\noindent {\bf References.}\setref{ABT1}
\medskip

\art ABT1 M. Abate, F. Bracci, F. Tovena: Index theorems for
holomorphic self-maps!   Ann. of Math.!  159 2004 819-864

\pre ABT2 M. Abate, F. Bracci, F. Tovena:  Index theorems for holomorphic maps and foliations!
Preprint, arXiv: math.CV/0601602!  2006 

\art An V. Ancona: Sur l'\'equivalence des voisinages des espaces analytiques contractibles!
Ann. Univ. Ferrara!  26  1980 165-172

\art At M. F. Atiyah: Complex analytic connections in fibre
bundles! Trans. Amer. Math. Soc.! 85 1957 181-207

\art BM A. Beauville, J.-Y. M\'erindol: Sections hyperplanes des surfaces $K3$! Duke Math. J.! 55
1987 873-878






\art CS C. Camacho, P. Sad: Invariant varieties through singularities of holomorphic
vector fields! Ann. of Math.! 115 1982 579-595


\book CM C. Camacho, H. Movasati: Neighborhoods of analytic varieties! Monograf\'\i as del
Instituto de Matem\'atica y Ciencias Afines,
35. Instituto de Matem\'atica y Ciencias Afines, IMCA, Lima, 2003. Cf. also arXiv:
math.CV/0208058v1


\art CMS C. Camacho, H. Movasati, P. Sad: Fibered neighborhoods of curves in surfaces! J.
Geom. Anal.! 13 2003 55-66


\coll CG M. Commichau, H. Grauert: Das formale Prinzip f\"ur kompakte komplexe
Untermannigfaltigkeiten mit 1-positivem Normalenb\"undel! Recent developments in several complex
variables, 1979! Ann. of Math. Studies 100, Princeton University Press, Princeton, 1981

\book Ei D. Eisenbud: Commutative Algebra with a view toward
algebraic geometry! Springer-Verlag, New York, 1994

\art G H. Grauert: \"Uber Modifikationen und exzeptionelle analytische Mengen! Math. Ann.! 146
1962 331-368

\art Gr P.A. Griffiths: The extension problem in complex analysis II; embeddings
with positive normal bundle! Amer. J. Math.! 88 1966 366-446

%
\book H R. Hartshorne: Algebraic geometry! Springer-Verlag, New
York, 1997

\art Hi A. Hirschowitz: On the convergence of formal equivalence between embeddings! Ann. of
Math.! 113 1981 501-514

\art K1 S. Kosarew: Ein allgemeines Kriterium f\"ur das formale Prinzip! J. Reine Angew. Math.!
388 1988 18-39

\coll K2 S. Kosarew: On some new results on the formal principle for
embeddings! Proceedings of the conference on algebraic geometry (Berlin, 1985)! 
Teubner-Texte Math., 92, Teubner, Leipzig, 1986, pp. 217--227


\book L H. B. Laufer: Normal two-dimensional singularities! Ann. Math. Studies 71, Princeton
University Press, Princeton, 1971





\art MR J. Morrow, H. Rossi: Submanifolds of ${\P}^n$ with
splitting normal bundle sequence are linear! Math. Ann.! 234 1978
253-261

\art MP M. Mustata, M. Popa: A new proof of a theorem of Van de
Ven! Bull. Math. Soc. Sc. Math. Roumanie! 39  1996 243-249

%

\art S M. Schlessinger: On rigid singularities! Rice Univ. Studies! 59 1973 147-162

\art St V. Steinbiss: Das formale Prinzip f\"ur reduzierte komplexe R\"aume mit einer
schwachen Positivit\"atsei\-gen\-schaft! Math. Ann.!  274  1986  485-502

%


%
\coll VdV A. Van de Ven: A property of algebraic varieties in
complex projective spaces! Colloque de g\'eom\'etrie
diff\'erentielle globale! Bruxelles, (1958) 151--152

\bye